\documentclass[amsfonts]{article}
\usepackage{amssymb}

\newtheorem{theorem}{Theorem}[section]
\newtheorem{corollary}[theorem]{Corollary}
\newtheorem{proposition}[theorem]{Proposition}
\newtheorem{lemma}[theorem]{Lemma}
\newtheorem{definition}[theorem]{Definition}
\newtheorem{remark}[theorem]{Remark}

\def\bR{\mathbb{R}}
\def\bZ{\mathbb{Z}}

\def\bH{\mathbb{H}}
\def\bQ{\mathbb{Q}}
\def\bM{\mathbb{M}}

\def\bD{\mathbb{D}}
\def\bL{\mathbb{L}}

\def\cB{\mathcal{B}}

\def\cD{\mathcal{D}}
\def\cE{\mathcal{E}}
\def\cF{\mathcal{F}}
\def\cG{\mathcal{G}}
\def\cH{\mathcal{H}}

\def\cS{\mathcal{S}}

\begin{document}

\title{$L_p$-Theory for the Stochastic Heat Equation
with Infinite-Dimensional Fractional Noise}

\author{R.M. Balan
\thanks{Research supported by a grant from the Natural Sciences and
Engineering Research Council of Canada.}  }

\date{April 9, 2009}

\maketitle

\begin{abstract}
\noindent In this article, we consider the stochastic heat equation
$du=(\Delta u+f(t,x))dt+ \sum_{k=1}^{\infty} g^{k}(t,x) \delta
\beta_t^k, t \in [0,T]$, with random coefficients $f$ and $g^k$,
driven by a sequence $(\beta^k)_k$ of i.i.d. fractional Brownian
motions of index $H>1/2$. Using the Malliavin calculus techniques
and a $p$-th moment maximal inequality for the infinite sum of
Skorohod integrals with respect to $(\beta^k)_k$, we prove that the
equation has a unique solution (in a Banach space of summability
exponent $p \geq 2$), and this solution is H\"older continuous in
both time and space.
\end{abstract}

\section{Introduction}

The study of stochastic partial differential equations driven by
colored noise has become an active area of research in the recent
years, which is viewed as an alternative (with an increased
potential for applications) to the classical theory of equations
perturbed by space-time white noise (see \cite{walsh86},
\cite{daprato-zabczyk92}, \cite{kallianpur-xiong95}, \cite{krylov99}
for fundamental developments -using different approaches- in the
white noise case.)

A Gaussian noise is said to be fractional in time, if its temporal
covariance structure coincides with that of a fractional Brownian
motion (fBm). Recall that a centered Gaussian process $(\beta_t)_{t
\in [0,T]}$ is a fBm of index $H \in (0,1)$ if
$R_{H}(t,s):=E(\beta_t \beta_s)=(t^{2H}+s^{2H}-|t-s|^{2H})/2$. The
case $H>1/2$ is referred as the ``regular'' case, whereas the case
$H=1/2$ corresponds to the Brownian motion. (The survey articles
\cite{nualart03} and \cite{hu05} offer more details on the fBm.)

Since the fBm is not a semimartingale, one cannot use the It\^o
calculus, which lies at the foundation of the study of equations
driven by white noise. Various methods exist in the literature to
circumvent this difficulty, based on the Skorokod integral (e.g.
\cite{alos-mazet-nualart01}, \cite{alos-nualart03},
\cite{carmona-coutin03}, \cite{decreusefond-ustunel99},
\cite{DHP00}), the pathwise generalized Stieltjes integrals (e.g.
\cite{zahle98}, \cite{nualart-vuillermont06},
\cite{sanzsole-vuillermont07}), or the ``rough paths'' analysis
(e.g.
 \cite{lyons98}, \cite{lyons-qian02}).

The present article is dedicated to the study of the stochastic heat
equation with (additive) infinite-dimensional fractional noise:
\begin{equation}
\label{heat} du(t,x)=(\Delta u(t,x)+f(t,x))dt+ \sum_{k=1}^{\infty}
g^{k}(t,x) \delta \beta_t^k, \quad t \in [0,T], x \in \bR^d,
\end{equation}
where $(\beta^k)_k$ is a sequence of i.i.d. fBm's of index $H>1/2$,
the solution is defined in the weak sense (using integration against
test functions $\phi \in C_0^{\infty}(\bR^d)$), and $\delta
\beta_t^k$ is a formal way of indicating that the stochastic
integrals (which are used for defining the solution) are interpreted
in the Skorohod sense.

Let $H_p^n(\bR^d)$ ($n \in \bR, p \geq 2$) be the Sobolev space of
all generalized functions on $\bR^d$ whose derivatives of order $k
\leq n$ lie in $L_p(\bR^d)$.
Our main result shows that for suitable initial condition $u_0$, and
Sobolev-space valued random processes $f=\{f(t, \cdot)\}_{t \in
[0,T]}$ and $g^k=\{g^k(t, \cdot)\}_{t \in [0,T]}, k \geq 1$,
equation (\ref{heat}) has a unique $H_p^n(\bR^d)$-valued solution
$u=\{u(t,\cdot)\}_{t \in [0,T]}$, and $u \in
C([0,T],H_p^{n-2}(\bR^d))$ a.s., such that
$$E\sup_{t \leq T}\|u(t, \cdot)\|_{H_p^{n-2}(\bR^d)}^p<\infty , \quad
E\int_0^T \|u(t,\cdot)\|_{H_p^n(\bR^d)}^pdt<\infty.$$

\noindent Moreover, $u$ belongs to the H\"older space
$C^{\alpha-1/p}([0,T], H_p^{n-2\beta})$, with probability $1$, for
any $1/2 \geq \beta>\alpha>1/p$. If in addition,
$\gamma:=n-2\beta-d/p>0$, $u$ is also $\gamma$-H\"older continuous
in space, since $H_p^{n-2\beta}(\bR^d) \subset C^{\gamma}(\bR^d)$.
These results provide generalizations to the fractional case of the
existing results for the heat equation driven by a sequence
$(w^k)_k$ of i.i.d. Brownian motions (see \cite{rozovskii90},
\cite{krylov96}, \cite{krylov99}).

We note that our result cannot be inferred from the results existing
in the literature for parabolic equations driven by Hilbert-space
valued fractional noise with trace-class covariance operator (e.g.
\cite{grecksch-anh99}, \cite{maslowski-nualart03}, \cite{TTV03}).
Nevertheless, we should mention the recent related investigations of
\cite{nualart-vuillermont06} and \cite{sanzsole-vuillermont07},
using fractional calculus techniques (as opposed to the Malliavin
calculus techniques used here), which establish the existence and
H\"older continuity (in time) of a variational/mild $L_2(D)$-valued
solution for a parabolic initial-boundary value problem with
multiplicative fractional noise, when $D \subset \bR^d$ is a bounded
open set.

Similarly to the Brownian motion case, at the origin of our
developments lie two basic tools: (1) a generalization of the
Littlewood-Paley inequality for Banach-space valued functions
(Theorem \ref{gener-Th1-1-Krylov}, Appendix); and (2) a suitable
$p$-th moment maximal inequality for the sum of Skorokod integrals
with respect to $(\beta^k)_k$ (Theorem \ref{max-inequality-p}):
\begin{eqnarray}
\nonumber \lefteqn{E\sup_{t \leq T} \left|\sum_{k=1}^{\infty}
\int_0^t u_s^k \delta \beta_s^k \right|^p \leq  C_{p,H,T} \left\{
E\left|\int_0^T \sum_{k=1}^{\infty}
|u_s^k|^{2}ds \right|^{p/2} + \right. } \\
\label{max} & & \left.E \left| \int_0^T \left[ \int_0^T
\left(\sum_{k=1}^{\infty}|D_{\theta}^{\beta^k}
u_s^k|^{2}\right)^{1/(2H)} d\theta \right]^{2H}ds \right|^{p/2}
\right\}.
\end{eqnarray}

Compared to the Burkholder-Davis-Gundy inequality (which was used in
the Brownian motion case), inequality (\ref{max}) contains an
additional term involving the Malliavin derivative $D^{\beta^k}u^k$
of the process $u^k$ with respect to $\beta^k$. It is because of
this extra term that our developments deviate significantly from the
white noise case, and we require that the multiplication coefficient
$g^k$ lie in a suitable space of Malliavin differentiable functions
with respect to $\beta^k$ (which in particular, implies that $g^k$
is measurable with respect to $\beta^k$).

This article is organized as follows. In Section 2, we give some
preliminaries on the Malliavin calculus for Hilbert-space valued
fractional processes, and we develop a maximal inequality for these
processes.
In Section 3, we convert the inequality obtained in Section 2 (which
speaks about the Skorohod integral with respect to a Hilbert-space
valued fractional process), into an inequality which speaks about
the sum of Skorohod intregrals with respect to a sequence
$(\beta^k)_k$ of i.i.d. fBm's. In Section 4, we introduce the
stochastic Banach spaces in which we are allowed to select the
coefficients $f$ and $(g^k)_k$. Section 5 is dedicated to the main
result, as well as the H\"older continuity of the solution.
The appendix contains the generalization of the Littlewood-Paley
inequality to Banach space valued functions.

\section{Malliavin Calculus for Fractional Processes}
\label{prelim-Malliavin}

In this section, we introduce the basic facts about the Malliavin
calculus with respect to (Hilbert-space valued) fractional
processes. We refer the reader to \cite{nualart98} and
\cite{nualart06} for a comprehensive account on this subject.
Throughout this work, we let $H \in (1/2,1)$ be fixed.

We begin by introducing some Banach spaces and Hilbert spaces of
deterministic functions, which are used for the Malliavin calculus
with respect to 
fractional processes.

If $V$ is an arbitrary Banach space, we let $\cE_V$ be the class of
all elementary functions $\phi:[0,T] \to V$ of the form
$\phi(t)=\sum_{i=1}^{m}1_{(t_{i-1},t_{i}]}(t)\varphi_i$ with $0 \leq
t_0 < \ldots <t_m \leq T$ and $\varphi_i \in V$. Let $|\cH_V|$ be
the space of all strongly measurable functions $\phi:[0,T] \to V$
with $\|\phi\|_{|\cH_V|}<\infty$, where
$$\|\phi\|_{|\cH_V|}^2:=\alpha_H \int_0^T \int_0^T \|\phi(t)\|_{V} \
\|\phi(s)\|_{V}  |t-s|^{2H-2} dtds, \quad \alpha_H=H (2H-1).$$

\noindent The space $\cE_{V}$ is dense in $|\cH_{V}|$ with respect
to the norm $\|\cdot \|_{|\cH_V|}$. It is known that there exists a
constant $b_H>0$ such that $\|\phi\|_{|\cH_V|} \leq b_H
\|\phi\|_{L_{1/H}([0,T];V)}$ for any $\phi \in L_{1/H}([0,T];V)$
(see e.g. relation (11) of \cite{alos-nualart03}).

In particular, if $V=\bR$, we denote $\cE_V=\cE$ and
$|\cH_V|=|\cH|$.

We let $|\cH| \otimes |\cH_V|$ be the space of all strongly
measurable functions $\phi:[0,T]^2 \to V$ with $\|\phi\|_{|\cH|
\otimes |\cH_V|}<\infty$, where
$$\|\phi\|_{|\cH| \otimes |\cH_V|}^{2}:=\alpha_{H}^2\int_{[0,T]^4}
\|\phi(t,\theta)\|_{V} \ \|\phi(s,\eta)\|_{V} \ |t-s|^{2H-2} \
|\theta-\eta|^{2H-2} d\theta d\eta dsdt.$$

If $V$ is a Hilbert space, we let $\cH_{V}$ be the completion of
$\cE_V$ with respect to the inner product $\langle \cdot, \cdot
\rangle_{\cH_{V}}$ defined by:
$$\langle \phi, \psi \rangle_{\cH_{V}}:=\alpha_H \int_0^T
\int_0^T \langle \phi(t), \psi(s) \rangle_{V}|t-s|^{2H-2} ds dt.$$
We have:
\begin{equation}
\label{ineq-norms-1} \|\phi\|_{\cH_{V}} \leq \|\phi \|_{|\cH_V|}
\leq b_{H}\| \phi \|_{L_{1/H}([0,T];V)} \leq b_H \|\phi
\|_{L_2([0,T];V)},
\end{equation}
and $L_{2}([0,T];V) \subset L_{1/H}([0,T];V)\subset |\cH_V| \subset
\cH_V$. In particular, if $V=\bR$, we denote $\cH_V=\cH$. The space
$\cH$ may
contain distributions of order $-(2H-1)$. 
Note that $\cH_V$ is isomorphic with $\cH \otimes V$, and the inner
products in the two spaces are the same.

We let $|\cH_V| \otimes |\cH_V|$ be the space of all strongly
measurable functions $\phi:[0,T]^2 \to V \otimes V$ with
$\|\phi\|_{|\cH_V| \otimes |\cH_V|}<\infty$, where
$$\|\phi\|_{|\cH_V| \otimes |\cH_V|}^{2}:=\alpha_{H}^2\int_{[0,T]^4}
\|\phi(t,\theta)\|_{V \otimes V} \ \|\phi(s,\eta)\|_{V \otimes V} \
|t-s|^{2H-2} \ |\theta-\eta|^{2H-2} d\theta d\eta dsdt,$$

\noindent and $\cH_V \otimes \cH_V$ be the completion of $\cE_V
\otimes \cE_V$ with respect to the inner product $\langle \cdot,
\cdot \rangle_{\cH_V \otimes \cH_V}$ defined by:
$$\langle \phi, \psi \rangle_{\cH_V \otimes \cH_V}:=\alpha_H^2
\int_{[0,T]^4} \langle \phi(t,\theta), \psi(s,\eta) \rangle_{V
\otimes V} |t-s|^{2H-2}|\theta-\eta|^{2H-2}d\theta d\eta ds dt.$$

We have: (see e.g Lemma 1, \cite{alos-nualart03} for the second
inequality below)
\begin{equation}
\label{ineq-norms-2} \|\phi\|_{\cH_V \otimes \cH_V} \leq \|\phi
\|_{|\cH_{V}| \otimes |\cH_V|} \leq b_{H}\| \phi
\|_{L_{1/H}([0,T]^2; V \otimes V)} \leq b_H \| \phi
\|_{L_{2}([0,T]^2; V \otimes V)},
\end{equation}
and $L_{2}([0,T]^2; V \otimes V) \subset L_{1/H}([0,T]^2; V \otimes
V) \subset |\cH_V| \otimes |\cH_V| \subset \cH_V \otimes \cH_V$.

\vspace{3mm}

We begin now to introduce the main ingredients of the Malliavin
calculus with respect to fractional processes.

Let $V$ be an arbitrary Hilbert space and $B=(B(\phi))_{\phi \in
\cH_{V}}$ be a centered Gaussian process, defined on a probability
space $(\Omega, \cF,P)$, with covariance:
\begin{equation}
\label{covariance-B} E(B(\phi)B(\psi))=\langle \phi, \psi
\rangle_{\cH_V}, \quad \forall \phi, \psi \in \cH_V.
\end{equation}
If we let $B_t(\varphi):=B(1_{[0,t]}\varphi)$ for any $\varphi \in
V,t \in [0,T]$, then
$$E(B_t(\varphi)B_s(\eta))=R_{H}(t,s)\langle \varphi, \eta
\rangle_V, \quad \forall \varphi, \eta \in V, s,t \in [0,T].$$ (In
particular, if $V=\bR$, then $\beta_t:=B(1_{[0,t]}), t \in [0,T]$ is
a fBm of index $H$.)

Let
$$\cS_{B}:=\{F=f(B(\phi_1),\ldots, B(\phi_n)); f \in
C_{b}^{\infty}(\bR^n), \phi_i \in \cH_V, n \geq 1\}$$ be the space
of all ``smooth cylindrical'' random variables,  where
$C_b^{\infty}(\bR^d)$ denotes the class of all bounded infinitely
differentiable functions on $\bR^n$, whose partial derivatives are
also bounded. Clearly $\cS_B \subset L_p(\Omega)$ for any $p \geq
1$.

The {\bf Malliavin derivative} of an element $F=f(B(\phi_1),\ldots,
B(\phi_n))\in \cS_{B}$, with respect to $B$, is defined by:
$$D^{B}F:=\sum_{i=1}^{n}\frac{\partial f}{\partial x_i}(B(\phi_1),\ldots,
B(\phi_n))\phi_i.$$ Note that $D^{B}F \in L_p(\Omega; \cH_V)$ for
any $p \geq 1$; by abuse of notation, we write $D^{B} F=(D_t^{B}
F)_{t \in [0,T]}$ even if $D_t^B F$ is not a function in $t$. We
endow $\cS_{B}$ with the norm:
$$\|F\|_{\bD_{B}^{1,p}}^{p}:=E|F|^p+E\|D^{\beta} F \|_{\cH_V}^{p},$$

\noindent and we let $\bD_{B}^{1,p}$ be the completion of $\cS_{B}$
with respect to this norm. The operator $D^{B}$ can be extended to
$\bD_{B}^{1,p}$. The adjoint
$$\delta^{B}: {\rm Dom} \ \delta^{B} \subset L_{2}(\Omega;
\cH_V) \to L_2(\Omega)$$ of the operator $D^{B}$, is called the {\bf
Skorohod integral} with respect to $B$. The operator $\delta^{B}$ is
uniquely defined by the following relation:
$$E(F \delta^{B}(U))=E \langle D^{B}F,U \rangle_{\cH}, \quad \forall
 F \in \bD^{1,2}_{B}.$$
Note that $E(\delta^{B}(U))=0$ for any $u \in {\rm Dom} \
\delta^{B}$. If $U \in {\rm Dom} \ \delta^{B}$, we use the notation
$U=(U_t)_{t \in [0,T]}$ and $\delta^{B}(U)=\int_0^T U_{s} \delta
B_{s}$.

If $V'$ is an arbitrary Hilbert space, we let
$$\cS_{B}(V'):=\{U=\sum_{j=1}^{m}F_j \phi_j; F_j \in \cS_{B}, \phi_j
\in V', m \geq 1\}$$

\noindent  be the class of all ``smooth  cylindrical'' $V'$-valued
random variables. Clearly $\cS_B(V') \subset L_p(\Omega;V')$ for any
$p \geq 1$.

The Malliavin derivative of an element $U=\sum_{j=1}^{m}F_j \phi_j
\in \cS_{B}(V')$ is defined by
$D^{B}U:=\sum_{j=1}^{m}(D^{B}F_j)\phi_j$. We have $D^BU \in
 L_p(\Omega; \cH_V \otimes V')$ for any $p \geq 1$. We endow $\cS_{B}(V')$ with the
norm:
$$\|U\|_{\bD_{B}^{1,p}(V')}^{p}:= E\|U\|_{V'}^p+E\|D^{B} U \|_{\cH_V
\otimes V'}^{p},$$
and let $\bD_{B}^{1,p}(V')$ be the completion of $\cS_{B}(V')$ with
respect to this norm. The operator $D^{B}$ can be extended to
$\bD_{B}^{1,p}(V')$.

In particular, if $V'=\cH_V$, then $\bD^{1,2}_{\beta}(\cH_V) \subset
{\rm Dom} \ \delta^{B}$. If $U \in \bD_{B}^{1,2}(\cH_V)$ then
$D^{B}U \in L_2(\Omega; \cH_V \otimes \cH_V)$; by abuse of notation,
we write $D^{B}U=(D_t^{B}U_s)_{s,t \in [0,T]}$.

The space $\bD_{B}^{1,2}(\cH_V)$ is viewed as a ``suitable'' class
of Skorohod integrands with respect to $B$. For any $U \in
\bD^{1,2}_{B}(\cH_V)$, we have:
\begin{eqnarray}
\nonumber E|\delta^{B}(U)|^2&=&E\|U\|_{\cH_V}^2+E(\langle D^{B}U,
(D^{B}U)^*\rangle_{\cH_V \otimes \cH_V)})\\
\label{L2-norm-Skorohod-integral} &\leq &
E\|U\|_{\cH_V}^2+E\|D^{B}U\|^2_{\cH_V \otimes
\cH_V}=\|U\|_{\bD_{B}^{1,2}(\cH_V)}^2,
\end{eqnarray}
where $(D^{B}U)^*$ is the adjoint of $D^{B}U$ in $\cH_V \otimes
\cH_V$.

The following result is a consequence of Meyer's inequalities.

\begin{proposition}[Proposition 2.4.4 of \cite{nualart98}]
\label{prop2-4-4-nualartLN} Let $p>1$ and $U \in
\bD_{B}^{1,p}(\cH_V)$. Then $U$ lies in the domain of $\delta^B$ in
$L_p(\Omega)$ and
$$E|\delta^B(U)|^p \leq C_{H,p} \{\|E(U)\|_{\cH_V}^p+ E\|D^B U\|_{\cH_V \otimes
\cH_V}^p\},$$ where $C_{H,p}$ is a constant depending on $H$ and
$p$.
\end{proposition}

As a consequence of Proposition \ref{prop2-4-4-nualartLN},
(\ref{ineq-norms-1}) and (\ref{ineq-norms-2}), we obtain:
\begin{equation}
\label{conseq-prop2-4-4-nualartLN} E|\delta^B(U)|^p \leq C_{H,p}
b_{H}\{\|E(U)\|_{L_{1/H}([0,T];V)}^p+
 E\|D^B U\|_{L_{1/H}([0,T]^2; V \otimes V)}^p\}.
\end{equation}

We denote by $\bD_{B}^{1,p}(|\cH_V|)$ the set of all elements $U \in
\bD^{1,p}_{B}(\cH_V)$, such that $U \in |\cH_V|$ a.s., $D^{B} U \in
|\cH_V| \otimes |\cH_V|$ a.s., and
$\|U\|_{\bD_{B}^{1,p}(|\cH_V|)}<\infty$,
where
\begin{eqnarray*}
\|U\|_{\bD_{B}^{1,p}(|\cH_V|)}^p&:=&E\|U\|_{|\cH_V|}^p+
E\|D^{B}U\|_{|\cH_V| \otimes |\cH_V|}^p.
\end{eqnarray*}

The following result generalizes Theorem 4 of \cite{alos-nualart03}
to the case of $V$-valued fractional processes.

\begin{theorem}
\label{theorem4-alos-nualart} Let $1/2<H<1$, $p>1/H$ and
$0<\varepsilon<H-1/p$. Then, there exists a constant $C$ depending
on $H,p,\varepsilon$ and $T$ such that
\begin{eqnarray}
\nonumber  E\sup_{t \leq T} \left|\int_0^t U_s \delta B_s \right|^p
& \leq & C \left\{ \left(\int_0^T \|E(U_s)\|_V^{1/(H-\varepsilon)}ds
\right)^{p(H-\varepsilon)} + \right.  \\
\label{theorem4-alos-nualart-relation} & & \left. E \left[ \int_0^T
\left( \int_0^T \|D_{\theta}^B U_s \|_{V \otimes V}^{1/H} d\theta
\right)^{\frac{H}{H-\varepsilon}}ds \right]^{p(H-\varepsilon)}
\right\}
\end{eqnarray}

\noindent for any process $U=(U_t)_{t \in [0,T]} \in
\bD_{B}^{1,p}(|\cH_V|)$ for which the right-hand side of
(\ref{theorem4-alos-nualart-relation}) is finite.
\end{theorem}

\noindent {\bf Proof:} The argument is similar to the one used in
the proof of Theorem 4 of \cite{alos-nualart03}. We include it for
the sake of completeness. Let $\alpha=1-1/p-\varepsilon$.

By writing $\int_0^t U_s \delta B_s=c_{\alpha} \int_0^t
(t-r)^{-\alpha} \left(\int_0^r U_s (r-s)^{\alpha-1} \delta B_s
\right) dr$, and using H\"older's inequality, we obtain:
$$E\sup_{t \leq T} \left|\int_0^t U_s \delta B_s \right|^p \leq c_{\alpha,p}
E \int_0^T  \left|\int_0^r U_{s} (r-s)^{\alpha-1}\delta
B_{s}\right|^{p}dr,$$ where $c_{\alpha,p}$ is a constant depending
on $\alpha$ and $p$. Using (\ref{conseq-prop2-4-4-nualartLN}), we
have:
\begin{eqnarray*}
\lefteqn{E\sup_{t \leq T} \left|\int_0^t U_s \delta B_s \right|^p
\leq  c_{\alpha,p,H} \left\{\int_0^T \left(\int_0^r
\|E(U_s)\|_{V}^{1/H}(r-s)^{(\alpha-1)/H}
ds \right)^{pH}dr \right. } \\
& & + \left. E \int_0^T \left( \int_0^r \int_0^T \|D_{\theta}^{B}U_s
\|_{V \otimes V}^{1/H}(r-s)^{(\alpha-1)/H} d\theta ds \right)^{pH}dr
\right\},
\end{eqnarray*}
where $c_{\alpha,p,H}$ is a constant which depends on $\alpha,p$
and$H$. The result follows by applying Hardy-Littlewood inequality
(p. 119 of \cite{stein70}). $\Box$

\vspace{2mm}

When $p \geq 2$, the previous theorem leads to the following result.

\begin{corollary}
\label{V-maximal-inequality} Let $1/2<H<1$ and $p \geq 2$ be
arbitrary. Then, there exists a constant $C$ depending on $H,p$ and
$T$ such that
\begin{eqnarray}
\nonumber  E\sup_{t \leq T} \left|\int_0^t U_s \delta B_s \right|^p
& \leq & C \left\{ \left(\int_0^T \|E(U_s)\|_V^{2}ds
\right)^{p/2} + \right.  \\
 \label{V-maximal-inequality-relation} & & \left. E
\left[ \int_0^T \left( \int_0^T \|D_{\theta}^B U_s \|_{V \otimes
V}^{1/H} d\theta \right)^{2H}ds \right]^{p/2} \right\}
\end{eqnarray}

\noindent for any process $U=(U_t)_{t \in [0,T]} \in
\bD_{B}^{1,p}(|\cH_V|)$ for which the right-hand side of
(\ref{V-maximal-inequality-relation}) is finite.
\end{corollary}

\noindent {\bf Proof:} The result follows by applying Theorem
\ref{theorem4-alos-nualart} with $\varepsilon< H-1/2$ 
and using the fact that $\|\phi \|_{L_{1/(H-\varepsilon)}([0,T])}
\leq C_T \|\phi \|_{L_2([0,T])}$ for any $\phi \in L_2([0,T])$.
$\Box$

\section{The Maximal Inequality}

The goal of this section is to translate the $p$-th moment maximal
inequality given by Corollary \ref{V-maximal-inequality} into a
similar inequality (in the $l_2$-norm) for a sequence $(u^k)_k$ of
Skorohod integrable processes, with respect to a sequence
$(\beta^k)_k$ of i.i.d. fBm's. The idea is
to recover a Gaussian process $B$ (as in Section 2) from
$(\beta^k)_k$, and to construct a Skorohod integrable process $U$
(with respect to $B$) from the sequence $(u^k)_k$, such that
$\delta^B(U 1_{[0,t]})=\sum_{k=1}^{\infty}
\delta^{\beta^k}(u^k1_{[0,t]})$ for all $t \in [0,T]$ a.s.

\vspace{3mm}

Let $\beta^k=(\beta_t^k)_{t \in [0,T]}, k \geq 1$ be a sequence of
i.i.d. fBm's of Hurst index $H>1/2$, defined on the same probability
space $(\Omega, \cF,P)$. Let $V$ be an arbitrary Hilbert space, and
$(e_k)_k$ a complete orthonormal system in $V$.

The first result shows that it is possible to construct a centered
Gaussian process $B$ with covariance (\ref{covariance-B}), from the
sequence $(\beta^k)_k$. This result is probably well-known; we state
it for the sake of completeness.

\begin{lemma}
\label{deterministic-integrands} Let $(\phi^k)_k \subset \cH$ be
such that $\sum_{k=1}^{\infty}\|\phi^k\|_{\cH}^2 <\infty$. Then:

a) $\varphi^{(N)}:=\sum_{k=1}^{N}\phi^k e_k \in \cH_V$ for all $N
\geq 1$, and there exists $\varphi:=\sum_{k=1}^{\infty}\phi^k e_k
\in \cH_V$ such that $\|\varphi^{(N)}-\varphi\|_{\cH_V} \to 0$ as $N
\to \infty$. We have:
\begin{equation}
\label{norm-H-HV} \|\varphi\|_{\cH_V}^2 =\sum_{k=1}^{\infty}
\|\phi^k\|_{\cH}^2;
\end{equation}

b) $B^{(N)}(\varphi):=\sum_{k=1}^{N}\beta^k(\phi^k)  \in
L_2(\Omega)$ for any $N \geq 1$, and there exists
$B(\varphi):=\sum_{k=1}^{\infty}\beta^k(\phi^k) \in L_2(\Omega)$
such that $E|B^{(N)}(\varphi)-B(\varphi)|^2 \to 0$ as $N \to
\infty$. The process $B=\{B(\varphi)\}_{\varphi \in \cH_V}$ is
Gaussian with mean zero and covariance (\ref{covariance-B}). In
particular, for any $t \in [0,T],\varphi \in V$, we have:
\begin{equation}
\label{delta-B-beta-k}
B_t(\varphi):=B(1_{[0,t]}\varphi)=\sum_{k=1}^{\infty} \langle
\varphi,e_k \rangle_V \beta_t^k \quad in \ L_2(\Omega).
\end{equation}
\end{lemma}

\noindent {\bf Proof:} a) The sequence $\{\varphi^{(N)}\}_N$ is
Cauchy in $\cH_V$, since
$(\varphi^{(N)}-\varphi^{(M)})(t)=\sum_{k=M+1}^{N}\phi^k(t)e_k$ for
any $N>M \geq 1$, and hence
\begin{eqnarray*}
\|\varphi^{(N)}-\varphi^{(M)}\|_{\cH_V}^2 & =&
\alpha_H  \sum_{k=M+1}^{N} \int_0^T \int_0^T
\phi^{k}(t)\phi^{k}(s)|t-s|^{2H-2}ds dt \\
&=& \sum_{k=M+1}^{N}\|\phi^k\|_{\cH}^2 \to 0, \ {\rm as} \ M,N \to
\infty.
\end{eqnarray*}

\noindent In particular,
$\|\varphi^{(N)}\|_{\cH_V}^2=\sum_{k=1}^{N}\|\phi^k\|_{\cH}^2$. By
letting $N \to \infty$, we obtain (\ref{norm-H-HV}).

b) The sequence $\{B^{(N)}(\varphi)\}_N$ is Cauchy in $L_2(\Omega)$,
since $B^{(N)}(\varphi)-B^{(M)}(\varphi)=\sum_{k=M+1}^N
\beta^k(\phi^k)$ for any $N>M \geq 1$, and hence
$$E|B^{(N)}(\varphi)-B^{(M)}(\varphi)|^2=\sum_{k=M+1}^{N}E|\beta^k(\phi^k)|^2=
\sum_{k=M+1}^{N}\|\phi^k\|_{\cH}^2 \to 0, \ {\rm as} \ M,N \to
\infty.$$

\noindent To prove (\ref{delta-B-beta-k}), note that
$1_{[0,t]}\varphi=\sum_{k=1}^{\infty} \phi^{k} e_k$, where
$\phi^k=1_{[0,t]}\langle \varphi,e_k \rangle_{V}$. It follows that
$B(1_{[0,t]}\varphi)=\sum_{k=1}^{\infty}
\beta^k(\phi^k)=\sum_{k=1}^{\infty} \langle \varphi,e_k \rangle_{V}
\beta_t^k$. $\Box$

\vspace{3mm}

We begin now to explore the relationship between the Malliavin
derivatives with respect to $(\beta^k)_k$ and the Malliavin
derivative with respect to $B$.

An immediate consequence of (\ref{delta-B-beta-k}) is that
$\beta_t^k=B(1_{[0,t]}e_k)$ for any $t \in [0,T]$, and hence
\begin{equation}
\label{beta-k-B-equality} \beta^k(\phi)=B(\phi e_k), \quad \forall
\phi \in \cH.
\end{equation}

Let $F=f(\beta^k(\phi_1), \ldots, \beta^k(\phi_n)) \in
\cS_{\beta^k}$ be arbitrary, with $f \in C_b^{\infty}(\bR^n)$ and
$\phi_i \in \cH$. Then $\varphi_i:=\phi_i e_k \in \cH_V$,
$F=f(B(\varphi_1), \ldots, B(\varphi_n)) \in \cS_B$, and
\begin{eqnarray*}
D_t^B F &=&\sum_{i=1}^{n}\frac{\partial f}{\partial
x_i}(B(\varphi_1), \ldots,
B(\varphi_n))\varphi_i=\left[\sum_{i=1}^{n}\frac{\partial
f}{\partial x_i}(\beta^k(\phi_1), \ldots, \beta^k(\phi_n))\phi_i
\right] e_k \\
&=&(D_{t}^{\beta^k}F)e_k.
\end{eqnarray*}

\noindent From here we conclude that $\cS_{\beta^k} \subset \cS_B$,
and for any $F \in \cS_{\beta^k}$,
$$\|D^B F\|_{\cH_V}=\|D^{\beta^k}F\|_{\cH}, \quad \|F
\|_{\bD_B^{1,p}}=\|F\|_{\bD_{\beta^k}^{1,p}}, \ \forall p \geq 1.$$

\noindent It follows that $\bD_{\beta^k}^{1,p} \subset
\bD_{B}^{1,p}$ for any $p \geq 1$, and
\begin{equation}
\label{D-beta-k-DB-agree} D^B F=(D^{\beta^k}F)e_k, \quad \mbox{for
any} \ F \in \bD_{\beta^k}^{1,2}.
\end{equation}

If $u=\sum_{j=1}^{m}F_j \phi_j \in \cS_{\beta^k}(\cH)$ is arbitrary,
with $F_j \in \cS_{\beta^k}$ and $\phi_j \in \cH$, then $ue_k \in
\bD_{B}^{1,2}(\cH_V)$ and
$$D^B(u e_k)=\sum_{j=1}^{m}(D^B F_j)\phi_j e_k=
\sum_{j=1}^{m}(D^{\beta^k} F_j)\phi_j e_k \otimes e_k=(D^{\beta^k}
u) e_k \otimes e_k.$$

\noindent In general, if $u \in \bD_{\beta^k}^{1,2}(\cH)$, then
$ue_k \in \bD_{B}^{1,2}(\cH_V)$ and
\begin{equation}
\label{Mal-deriv-beta-B} D^B(u e_k)=(D^{\beta^k} u) e_k \otimes e_k.
\end{equation}

\noindent Moreover, we have the following result:

\begin{lemma}
\label{Mal-deriv-beta-B-sum} If $u^k \in \bD_{\beta^k}^{1,2}(\cH)$,
then $\sum_{k=1}^{N}u^ke_k \in \bD_{B}^{1,2}(\cH_V)$, $D^B
(\sum_{k=1}^{N}u^ke_k)=\sum_{k=1}^{N}(D^{\beta^k} u^k)e_k \otimes
e_k$, and $\|\sum_{k=1}^{N}u^k
e_k\|_{\bD_{B}^{1,2}(\cH_V)}^2=\sum_{k=1}^{N}
\|u^k\|_{\bD_{\beta^k}^{1,2}(\cH)}^{2}$.
\end{lemma}

\noindent {\bf Proof:} The result follows from the definitions of
the norms in $\bD_{B}^{1,2}(\cH_V)$, respectively
$\bD_{\beta^k}^{1,2}(\cH)$, and the following two identities:
\begin{eqnarray*}
\|\sum_{k=1}^{N}u^k e_k\|_{\cH_V}^2 & = & \alpha_H \int_0^T \int_0^T
\langle \ \sum_{k=1}^{N}u_t^k e_k, \sum_{l=1}^{N}u_s^l e_l \
\rangle_V
|t-s|^{2H-2}ds dt \\
&=& \alpha_H \sum_{k,l=1}^{N}  \int_0^T \int_0^T u_t^k u_s^l \langle
e_k,e_l \rangle_{V} |t-s|^{2H-2}ds dt \\
&=& \sum_{k=1}^{N} \|u^k \|_{\cH}^2 \\
\|D^B(\sum_{k=1}^{N}u^k e_k)\|_{\cH_V \otimes \cH_V}^2 & = &
\alpha_H^2 \int_{[0,T]^4}\langle \
\sum_{k=1}^{N}D_{\theta}^B(u_{t}^k
e_k), \sum_{l=1}^{N}D_{\eta}^B(u_{s}^l e_l) \ \rangle_{V \otimes V} \\
& & |t-s|^{2H-2}|\theta-\eta|^{2H-2} d\theta d\eta ds dt \\
& = & \alpha_H^2 \sum_{k,l=1}^{N}  \int_{[0,T]^4}
(D_{\theta}^{\beta^k} u_{t}^k) (D_{\eta}^{\beta^k}u_{s}^l) \langle
e_k \otimes e_k,e_l
\otimes e_l \rangle_{V \otimes V} \\
& & |t-s|^{2H-2}
|\theta-\eta|^{2H-2} d\theta d\eta ds dt \\
&=& \sum_{k=1}^{N} \|D^{\beta^k}u^k \|_{\cH \otimes \cH}^2,
\end{eqnarray*}

\noindent where we used (\ref{Mal-deriv-beta-B}) for the second-last
equality above. $\Box$

\vspace{3mm}

We need an auxiliary result.

\begin{lemma}
\label{elementary-lemma} Let $X$ be a normed space and $y_N,
x_{N,n},x_n, x \in X$ be such that: $\lim_{n \to \infty}\sup_{N \geq
1}\|y_N-x_{N,n}\|=0$, $\lim_{N \to \infty}\|x_{N,n}-x_n \|=0$ for
all $n$, and $\lim_{n \to \infty}\|x_n-x\|=0$. Then $\lim_{N \to
\infty}\|y_N-x\|=0$.
\end{lemma}

\noindent {\bf Proof:} We use $\|y_N-x\| \leq
\|y_N-x_{N,n}\|+\|x_{N,n}-x_n\|+\|x_n-x\|$. $\Box$

\vspace{3mm}

The previous observations allow us to extend Lemma
\ref{deterministic-integrands} to the case of random integrands.

\begin{theorem}
\label{connection} Let $u^k \in \bD_{\beta^k}^{1,2}(\cH)$ for all $k
\geq 1$, such that
\begin{equation}
\label{cond-on-uk} \sum_{k=1}^{\infty} \|u^k\|_{
\bD_{\beta^k}^{1,2}(\cH)}^{2}<\infty.
\end{equation}
 Then:

a) $U^{(N)}:=\sum_{k=1}^{N}u^k e_k \in \bD_{B}^{1,2}(\cH_V)$ for any
$N \geq 1$, and there exists $U:=\sum_{k=1}^{\infty} u^k e_k \in
\bD_{B}^{1,2}(\cH_V)$ such that
$\|U^{(N)}-U\|_{\bD_{B}^{1,2}(\cH_V)} \to 0$ as $N \to \infty$. We
have: $D^B U=\sum_{k=1}^{\infty} (D^{\beta^k} u^k) e_k \otimes e_k$
and
\begin{equation}
\label{norm-H-HV-for-U}
\|U\|_{\bD_{B}^{1,2}(\cH_V)}^2=\sum_{k=1}^{\infty}
\|u^k\|_{\bD_{\beta^k}^{1,2}(\cH)}^2.
\end{equation}

b) the sequence $W^{(N)}:=\sum_{k=1}^{N}\delta^{\beta^k}(u^k), N
\geq 1$ has a limit in $L_2(\Omega)$, which coincides with
$\delta^B(U)$. We write
\begin{equation}
\label{Skorohod-integrals-agree} \delta^B(U)=\sum_{k=1}^{\infty}
\delta^{\beta^k}(u^k) \quad in \ L_2(\Omega).
\end{equation}
\end{theorem}

\noindent {\bf Proof:} a) By Lemma \ref{Mal-deriv-beta-B-sum},
$\{U^{(N)}\}_N$ is a Cauchy sequence in $\bD_{B}^{1,2}(\cH_V)$,
since:
$$\|U^{(N)}-U^{(M)}\|_{\bD_{B}^{1,2}(\cH_V)}^{2}=\sum_{k=M+1}^{N}
\|u^{k}\|_{\bD_{\beta^k}^{1,2}(\cH)}^{2} \to 0, \ {\rm as} \ M,N \to
\infty.$$  Hence, $U:=\lim_{N \to \infty} U^{(N)}$ exists in
$\bD_{B}^{1,2}(\cH_V)$, and $D^{B}U=\lim_{N \to \infty}D^{B}U^{(N)}$
in $L_2(\Omega; \cH_V \otimes \cH_V)$. Also,
$\|U^{(N)}\|_{\bD_{B}^{1,2}(\cH_V)}^{2}=\sum_{k=1}^{N}
\|u^{k}\|_{\bD_{\beta^k}^{1,2}(\cH)}^{2}$, and relation
(\ref{norm-H-HV-for-U}) follows by letting $N \to \infty$.

b) By inequality (\ref{L2-norm-Skorohod-integral}) (applied for
$V=\bR$ and $B=\beta^k$), we have:
$$\sum_{k=M+1}^{N}E|\delta^{\beta^k}(u^k)|^2 \leq
\sum_{k=M+1}^{N}\|u^k\|_{\bD_{\beta^k}^{1,2}(\cH)}^2 \to 0, \ {\rm
as} \ M,N \to \infty,$$ i.e. the sequence $\{W^{(N)}\}_N$ is Cauchy
in $L_2(\Omega)$. We let $W$ be the limit of $\{W^{(N)}\}_N$ in
$L_2(\Omega)$. We now prove that $W=\delta^B(U)$ (in $L_2(\Omega)$).

{\em Step 1.} Suppose that $u^k \in \cS_{\beta^k}(\cH)$ for all $k$,
i.e. $u^k=\sum_{j=1}^{m_k}F_j^k \phi_j^k$ for some $F_j^k \in
\cS_{\beta^k}$ and $\phi_j^k \in \cH$. Since $U^{(N)} \to U$ in
$\bD_{B}^{1,2}(\cH_V)$, $\delta^B(U^{(N)}) \to \delta^B(U)$ in
$L_2(\Omega)$. On the other hand
$\sum_{k=1}^{N}\delta^{\beta^k}(u^k) \to W$ in $L_2(\Omega)$. Hence,
it suffices to prove that:
\begin{equation}
\label{delta-Bbeta-agree-simple}
\delta^B(U^{(N)})=\sum_{k=1}^{N}\delta^{\beta^k}(u^k).
\end{equation}

\noindent Note that $U^{(N)}=\sum_{k=1}^{N}\sum_{k=1}^{m_k}F_j^k
\phi_j^k e_k \in \cS_B(\cH_V)$, since $F_j^k \in \cS_{\beta^k}
\subset S_{B}$ and $\phi_j^k e_k \in \cH_V$. Relation
(\ref{delta-Bbeta-agree-simple}) follows from
(\ref{beta-k-B-equality}) and (\ref{D-beta-k-DB-agree}), since:
\begin{eqnarray*}
\delta^B (U^{(N)})&=&\sum_{k=1}^{N}\sum_{j=1}^{m_k}F_j^k B(\phi_j^k
e_k)-\sum_{k=1}^{N}\sum_{j=1}^{m_k}\langle D^{B}F_{j}^k, \phi_j^ke_k
\rangle_{\cH_V} \\
\sum_{k=1}^{N}\delta^{\beta^k}(u^k) &=&
\sum_{k=1}^{N}\sum_{j=1}^{m_k}F_j^k
\beta^k(\phi_j^k)-\sum_{k=1}^{N}\sum_{j=1}^{m_k}\langle
D^{\beta^k}F_{j}^k, \phi_j^k \rangle_{\cH}.
\end{eqnarray*}
(We used relation (1.9) of \cite{nualart98}, for the equalities
above.)

{\em Step 2.} Suppose that $u^k \in \bD_{\beta^k}^{1,2}(\cH)$ for
all $k$. For any $\varepsilon>0$, there exists $u_{\varepsilon}^k
\in \cS_{\beta^k}(\cH)$ such that
$\|u_{\varepsilon}^k-u^k\|_{\bD_{\beta^k}^{1,2}(\cH)}<\varepsilon/2^k$;
hence
$\sum_{k=1}^{\infty}\|u_{\varepsilon}^k\|_{\bD_{\beta^k}^{1,2}(\cH)}^{2}<\infty$.
By part a), $U_{\varepsilon}:=\sum_{k=1}^{\infty} u_{\varepsilon}^k
e_k \in \bD_{B}^{1,2}(\cH_V)$ and
$\|U_{\varepsilon}-U\|_{\bD_{B}^{1,2}(\cH_V)}^2 =\sum_{k=1}^{\infty}
\|u_{\varepsilon}^k-u^k\|_{\bD_{\beta^k}^{1,2}(\cH)}^2 \leq
\varepsilon^2$. Taking $\varepsilon=1/n$, we conclude that for any
$k$, there exists a sequence $(u_n^k)_n \subset \cS_{\beta^k}(\cH)$,
such that $U_{n}:=\sum_{k=1}^{\infty} u_{n}^k e_k \in
\bD_{B}^{1,2}(\cH_V)$ and
$$\|U_{n}-U\|_{\bD_{B}^{1,2}(\cH_V)}^2=\sum_{k=1}^{\infty}
\|u_{n}^k-u^k\|_{\bD_{\beta^k}^{1,2}(\cH)} ^2 \to 0, \quad \mbox{as}
\ n \to \infty.$$

We now invoke Lemma \ref{elementary-lemma}, with $X=L_2(\Omega)$,
and
$$y_N=W^{(N)}=\sum_{k=1}^{N} \delta^{\beta^k}(u^k), \quad
x_{N,n}=\sum_{k=1}^{N}\delta^{\beta^k}(u_n^k), \quad
x_n=\delta^B(U_n), \quad x=\delta^B(U).$$ The hypothesis of the
lemma are verified, since $\lim_{N \to
\infty}\|x_{N,n}-x_n\|_{L_2(\Omega)}=0$ for all $n$ (by {\em Step
1}), $\lim_{n \to \infty} \|x_n-x\|_{L_2(\Omega)}=0$ (since $U_n \to
U$ in $\bD_{B}^{1,2}(\cH_V)$),
\begin{eqnarray*}
\|y_{N}-x_{N,n}\|_{L_2(\Omega)}^{2}&=& E
\left|\sum_{k=1}^{N}\delta^{\beta^k}(u^k-u_n^k)\right|^2 =
\sum_{k=1}^{N} E |\delta^{\beta^k}(u^k-u_n^k)|^2 \\
& \leq & \sum_{k=1}^{N} \|u^k-u_n^k\|_{\bD_{\beta^k}^{1,2}(\cH)}^2,
\end{eqnarray*}
and hence $\sup_{N \geq 1} \|y_{N}-x_{N,n}\|_{L_2(\Omega)}^{2} \leq
\sum_{k=1}^{\infty} \|u^k-u_n^k\|_{\bD_{\beta^k}^{1,2}(\cH)}^2 \to
0, \ {\rm as} \ n \to \infty$. We conclude that $\lim_{N \to
\infty}\|y_N-x\|_{L_2(\Omega)}=0$, i.e. $W=\delta^{B}(U)$. $\Box$

\vspace{3mm}

In the case $p=2$, we have the following preliminary result.

\begin{theorem}
\label{max-inequality-p-equal-2} There exists a constant $C$
depending on $H$ and $T$ such that
\begin{eqnarray}
\nonumber  E\sup_{t \leq T} \left|\sum_{k=1}^{\infty} \int_0^t u_s^k
\delta \beta_s^k \right|^2 & \leq & C \left\{ \sum_{k=1}^{\infty}E
\int_0^T |u_s^{k}|^{2}ds
 + \right.  \\
 \label{max-inequality-p-equal-2-relation} & & \left. \sum_{k=1}^{\infty} E
\int_0^T \left( \int_0^T |D_{\theta}^{\beta^k} u_s^k |^{1/H} d\theta
\right)^{2H}ds  \right\}
\end{eqnarray}

\noindent for any process $u=(u^k)_k$ for which $u^k \in
\bD_{\beta^k}^{1,2}(|\cH|)$ for all $k \geq 1$, and the right-hand
side of (\ref{max-inequality-p-equal-2-relation}) is finite.
\end{theorem}

\noindent {\bf Proof:} Let $0<\varepsilon<H-1/2$ and
$\alpha=1/2-\varepsilon$. As in the proof of Theorem 4,
\cite{alos-nualart03}, one can show that $$\sup_{t \leq T}
\left|\sum_{k=1}^{\infty}\int_0^t u_s^k \delta \beta_s^k \right|^2
\leq c_{\alpha}'\int_0^T \left| \sum_{k=1}^{\infty}\int_0^r u_s^k
(r-s)^{\alpha-1}\delta \beta_s^k\right|^2 dr.$$

\noindent Since the random variables $X_k=\int_0^r
u_s^{k}(r-s)^{\alpha-1} \delta \beta_s^k, k \geq 1$ are independent
with zero mean, $E(\sum_{k=1}^{n}X_k)^2=\sum_{k=1}^{n}E(X_k^2)$ for
all $n$. By the Fatou's lemma, $$ E\left|
\sum_{k=1}^{\infty}\int_0^r u_s^k (r-s)^{\alpha-1}\delta
\beta_s^k\right|^2 \leq \sum_{k=1}^{\infty} E \left|\int_0^r u_s^k
(r-s)^{\alpha-1}\delta \beta_s^k\right|^2.$$

Using (\ref{conseq-prop2-4-4-nualartLN}) and H\"older's inequality
we get:
\begin{eqnarray*}
\lefteqn{E \sup_{t \leq T}\left|\sum_{k=1}^{\infty}\int_0^t u_s^k
\delta \beta_s^k \right|^2  \leq  c_{\alpha}' \sum_{k=1}^{\infty}
\int_0^T E \left| \int_0^r u_s^k(r-s)^{\alpha-1}\delta \beta_s^k
\right|^{2}dr } \\
& &  \leq c_{\alpha,H} \left\{\sum_{k=1}^{\infty}\int_0^T E \left(
\int_0^r |u_s^k|^{1/H} (r-s)^{(\alpha-1)/H}ds \right)^{2H}dr+
\right. \\
& & \left. \sum_{k=1}^{\infty}  \int_0^T  E \left( \int_0^r\int_0^T
|D_{\theta}^{\beta^k}u_s^k|^{1/H} d\theta (r-s)^{(\alpha-1)/H}
ds\right)^{2H}dr \right\}  \\
& & \leq c_{\alpha,H} \left\{ \sum_{k=1}^{\infty} \int_0^T
r^{2(\alpha-1)+2H-1} E \int_0^r |u_s^k|^2ds dr + \right.
\\
& & \left. \sum_{k=1}^{\infty} \int_0^T r^{2(\alpha-1)+2H-1}  E
\int_0^r \left( \int_0^T |D_{\theta}^{\beta^k} u_s^k|^{1/H}d\theta
\right)^{2H}ds dr \right\}.
\end{eqnarray*}
$\Box$

 Let $l_2$ be the set of sequences $a=(a^k)_k, a^k \in \bR$ with
$|a|_{l_2}^2:=\sum_{k=1}^{\infty} |a^k|^2<\infty$. If $u=(u^k)_k$ is
such that $u^k \in \bD_{\beta^k}^{1,2}(\cH)$ for all $k \geq 1$, we
denote $Du:=(D^{\beta^k}u^k)_k$.

\vspace{3mm}

The next theorem is the main result of this section. Its proof is
based on Corollary \ref{V-maximal-inequality}, the connection
between the Skorohod integrals with respect to $(\beta^k)_k$ and the
Skorohod integral with respect to $B$ (given by Theorem
\ref{connection}), and Theorem \ref{max-inequality-p-equal-2}.

\begin{theorem}
\label{max-inequality-p} Let $1/2<H<1$ and $p \geq 2$. Then, there
exists a constant $C$ depending on $H,p$ and $T$ such that
\begin{eqnarray}
\nonumber  E\sup_{t \leq T} \left|\sum_{k=1}^{\infty} \int_0^t u_s^k
\delta \beta_s^k \right|^p & \leq & C \left\{ E\left(\int_0^T
|u_s|_{l_2}^{2}ds
\right)^{p/2} + \right.  \\
 \label{max-inequality-p-relation} & & \left. E
\left[ \int_0^T \left( \int_0^T |D_{\theta} u_s|_{l_2}^{1/H} d\theta
\right)^{2H}ds \right]^{p/2} \right\}
\end{eqnarray}

\noindent for any process $u=(u^k)_k$ for which $u^k \in
\bD_{\beta^k}^{1,p}(|\cH|)$ for all $k \geq 1$, and the right-hand
side of (\ref{max-inequality-p-relation}) is finite.
\end{theorem}

\noindent {\bf Proof:} Let $u=(u^k)_k$ be such that $u^k \in
\bD_{\beta^k}^{1,p}(|\cH|)$ for all $k \geq 1$, and the right-hand
side of (\ref{max-inequality-p-relation}) is finite. Since $p \geq
2$, $|EX|^{p/2} \leq E|X|^{p/2}$, for any $X \in L_{p/2}(\Omega)$,
and hence, $E \int_0^T|u_s|_{l_2}^{2}ds<\infty$ and $E \int_0^T
\left( \int_0^T |D_{\theta} u_s|_{l_2}^{1/H} d\theta \right)^{2H}ds
<\infty$. By Minkowski's inequality, $\sum_{k=1}^{\infty}  E
\int_0^T \left( \int_0^T |D_{\theta}^{\beta^k} u_s^k|^{1/H} d\theta
\right)^{2H} ds <\infty$. From here we conclude that relation
(\ref{cond-on-uk}) holds, since:
\begin{eqnarray*}
\sum_{k=1}^{\infty}\|u^{k}\|_{\bD_{\beta^k}^{1,2}(\cH)}^2
&=&\sum_{k=1}^{\infty}E\|u^k\|_{\cH}^2+\sum_{k=1}^{\infty}
E\|D^{\beta^k}u^k\|_{\cH \otimes
\cH}^2 \\
& \leq &
\sum_{k=1}^{\infty}E\|u^k\|_{L_2([0,T])}^2+\sum_{k=1}^{\infty}
E\|D^{\beta^k}u^k\|_{L_{2} ([0,T]; L_{1/H}([0,T]))}^2 <\infty.
\end{eqnarray*}

\noindent By Theorem \ref{connection}.(a), there exists
$U:=\sum_{k=1}^{\infty} u^k e_k \in \bD_{B}^{1,2}(\cH_V)$ and
$D^{B}U=\sum_{k=1}^{\infty} (D^{\beta^k}u^k)e_k \otimes e_k$.
Similarly, $U1_{[0,t]}=\sum_{k=1}^{\infty} u^k 1_{[0,t]} e_k \in
\bD_{B}^{1,2}(\cH_V)$ for any $t \in [0,T]$.

For any $t \in [0,T]$, let $$X_t:=\sum_{k=1}^{\infty}
\delta^{\beta^k}(u^k 1_{[0,t]}) \quad \mbox{and} \quad
Y_t:=\delta^B(U1_{[0,t]}).$$ Using the same argument as in Theorem 5
of \cite{alos-nualart03}, one can prove that $Y=(Y_t)_{t \in [0,T]}$
has an a.s. continuous modification. We work with this modification.

Also, for each $N \geq 1$, the process $X^{(N)}=(X_t^{(N)})_{t \in
[0,T]}$, defined by $X_t^{(N)}:=\sum_{k=1}^{N}\delta^{\beta^k}(u^k
1_{[0,t]}), t \in [0,T]$, has an a.s. continuous modification.

By Chebyshev's inequality, Theorem \ref{max-inequality-p-equal-2},
and (\ref{cond-on-uk}), the sequence $(X^{(N)})_N$ converges in
probability to $X$, in the sup-norm metric, since for any
$\varepsilon>0$,
\begin{eqnarray}
\label{convergence-in-sup-metric} \lefteqn{P(\sup_{t \leq T}
|X_t^{(N)}-X_t|>\varepsilon) \leq \frac{1}{\varepsilon^2} E \sup_{t
\leq T}\left|\sum_{k=N+1}^{\infty}\int_0^t u_s^k \delta
\beta_s^k\right|^2 \leq } \\
\nonumber & &  \frac{C}{\varepsilon^2} \left\{
\sum_{k=N+1}^{\infty}E \int_0^T |u_s^{k}|^{2}ds  +
\sum_{k=N+1}^{\infty} E \int_0^T \left( \int_0^T
|D_{\theta}^{\beta^k} u_s^k |^{1/H} d\theta \right)^{2H}ds \right\}
\to 0,
\end{eqnarray}
as $N \to \infty$.
Therefore, $X$ has an a.s. continuous modification. We work with
this modification.

From Theorem \ref{connection}.(b), we know that $Y_t=X_t$ a.s., for
any $t \in [0,T]$. Since both $Y$ and $X$ are a.s. continuous, it
follows that $Y_t=X_t$ for all $t \in [0,T]$ a.s. In particular,
$E\sup_{t \leq T}|Y_t|^p =E\sup_{t \leq T}|X_t|^p$, i.e.
\begin{equation}
\label{sup-delta-Bbeta-agree} E \sup_{t \leq T} \left|\int_0^t U_s
\delta B_s \right|^p=E\sup_{t \leq T} \left|\sum_{k=1}^{\infty}
\int_0^t u_s^k \delta \beta_s^k \right|^p.
\end{equation}

We now invoke Corollary \ref{V-maximal-inequality}. Note that
$E(U_s)=\sum_{k=1}^{\infty} E(u_s^k)e_k$. Hence $\|E(U_s)\|_V^2=
\sum_{k=1}^{\infty}|E(u_s^k)|^2 \leq \sum_{k=1}^{\infty}
E|u_s^k|^2=E|u_s|_{l_2}^2$ for any $s \in [0,T]$, and
\begin{equation}
\label{first-term-U-u} \left( \int_0^T \|E(U_s)\|_{V}^2
ds\right)^{p/2} \leq \left( E \int_0^T |u_s|_{l_2}^2 ds\right)^{p/2}
\leq E \left(\int_0^T |u_s|_{l_2}^2ds \right)^{p/2}.
\end{equation}

\noindent Note also that $D_{\theta}^B U_s=\sum_{k=1}^{\infty}
(D_{\theta}^{\beta^k}u_s^k) e_k \otimes e_k$, and hence,
\begin{equation}
\label{second-term-U-u} \|D_{\theta}^B U_s \|_{V \otimes
V}^2=\sum_{k=1}^{\infty}
|D_{\theta}^{\beta^k}u_s^k|^2=|D_{\theta}u|_{l_2}^2.
\end{equation}

\noindent Relation (\ref{max-inequality-p-relation}) becomes a
consequence of (\ref{V-maximal-inequality-relation}), combined with
(\ref{sup-delta-Bbeta-agree}), (\ref{first-term-U-u}) and
(\ref{second-term-U-u}). $\Box$

The following result is an immediate consequence of Theorem
\ref{max-inequality-p}.

\begin{corollary}
\label{max-inequality-p-cor} Let $1/2<H<1$ and $p \geq 2$. Then,
there exists a constant $C$ depending on $H,p$ and $T$ such that
\begin{eqnarray}
\nonumber \lefteqn{E\sup_{t \leq T} \left|\sum_{k=1}^{\infty}
\int_0^t u_s^k \delta \beta_s^k \right|^p \leq  C \left\{ E \int_0^T
|u_s|_{l_2}^{p}ds + \right. } \\
 \label{max-inequality-p-cor-relation} & & \left. E
\int_0^T \left( \int_0^T |D_{\theta} u_s|_{l_2}^{1/H} d\theta
\right)^{pH}ds  \right\}:=C \|u\|_{{\bL}_{H}^{1,p}(l_2)}^p
\end{eqnarray}

\noindent for any process $u=(u^k)_k$ for which $u^k \in
\bD_{\beta^k}^{1,p}(|\cH|)$ for all $k \geq 1$, and the right-hand
side of (\ref{max-inequality-p-cor-relation}) is finite.
\end{corollary}

\section{Stochastic Banach Spaces}

In this section, we introduce some Banach spaces of stochastic
integrands for the sequence of Skorohod integrals with respect to
$(\beta^k)_k$, which are suitable for our analysis. To ease the
exposition, we first treat the case of a single fBm (subsection
\ref{single-fBm}), and then the case of a sequence of i.i.d. fBm's
(subsection \ref{sequence-fBms}).

\subsection{The case of a single fBm}
\label{single-fBm}

We begin by recalling some basic facts about fractional Sobolev
spaces, using the notation in \cite{krylov99}. We let
$C_0^{\infty}=C_0^{\infty}(\bR^d)$ be the space of infinitely
differentiable functions on $\bR^d$, with compact support, and
$\cD=\cD(\bR^d)$ be the space of real-valued Schwartz distributions
on $C_0^{\infty}$. For $p \geq 1$, we denote by $L_p=L_p(\bR^d)$ the
set of all measurable functions $u:\bR^d\to \bR$ such that
$\|u\|_{L_p}^p:=\int_{\bR^d}|u(x)|^p dx<\infty$.

For any $p >1$ and $n \in \bR$, we let $H_p^n=H_p^n(\bR^d):=\{u \in
\cD; (1-\Delta)^{n/2}u \in L_p\}$ be the fractional Sobolev space,
with the norm $\|u\|_{H_p^n}:=\|(1-\Delta)^{n/2}u\|_{L_p}$. For any
$u \in H_p^n$ and $\phi \in C_0^{\infty}$, we define
$$(u,\phi):=\int_{\bR^d}[(1-\Delta)^{n/2}u](x) \cdot
[(1-\Delta)^{-n/2}\phi](x) dx.$$

\noindent By H\"older's inequality, for any $u \in H_p^n$ and $\phi
\in C_0^{\infty}$, we have:
\begin{equation}
\label{norm-of-g-Hpn} |(u,\phi)|^2 \leq N \|u \|_{H_p^n}^2,
\end{equation}
where
$N= \|(1-\Delta)^{-n/2}\phi \|_{L_{p/(p-1)}}^2$ is a constant
depending on $n,p$ and $\phi$.

\vspace{2mm}

Let $\beta=(\beta_t)_{t \in [0,T]}$ be a fBm of index $H>1/2$,
defined on a probability space $(\Omega,\cF,P)$. We introduce the
following spaces of Banach-space valued integrands for the Skorohod
integral with respect to $\beta$.

\begin{definition} Let $V$ be an arbitrary Banach space and $p >1$.

a) We denote by $\bD_{\beta}^{1,p}(|\cH_{V}|)$ the set of all
elements $g \in \bD_{\beta}^{1,p}(\cH_{V})$ such that $g \in
|\cH_{V}| $ a.s., $D^{\beta} g \in |\cH| \otimes |\cH_{V}|$ a.s.,
and $\|g\|_{\bD_{\beta}^{1,p}(|\cH_{V}|)}<\infty$, where
$$\|g\|_{\bD_{\beta}^{1,p}(|\cH_{V}|)}^p :=
E\|g\|_{|\cH_{V}|}^p+E\|D^{\beta}g\|_{|\cH| \otimes
|\cH_{V}|}^{p}.$$

b) We denote by $\bL_{H,\beta}^{1,p}(V)$ the set of all elements $g
\in \bD_{\beta}^{1,p}(|\cH_V|)$ such that $\|g
\|_{\bL_{H,\beta}^{1,p}(V)}<\infty$, where
$$\|g\|_{\bL_{H,\beta}^{1,p}(V)}^p := E\int_0^T \|g_s\|_{V}^pds+E
\int_0^T \left( \int_0^T \|D_{t}^{\beta}g_s\|_{V}^{1/H} dt
\right)^{pH}ds.$$

c) We denote by $\widetilde{\bL}_{H,\beta}^{\ 1,p}(V)$ the
completion of $\cS_{\beta}(\cE_{V})$ in
$\bD_{\beta}^{1,p}(|\cH_{V}|)$, with respect to the norm
$\|\cdot\|_{{\bL}_{H,\beta}^{1,p}(V)}$.
\end{definition}

\noindent Using (\ref{ineq-norms-1}) and (\ref{ineq-norms-2}), one
can prove that:
\begin{equation}
\label{ineq-norms-L-D-U} \|g\|_{\bD_{\beta}^{1,p}(|\cH_V|)} \leq b_H
\|g\|_{\bL_{H,\beta}^{1,p}(V)}, \quad \forall u \in
\bL_{H,\beta}^{1,p}(V).
\end{equation}

\begin{remark}
{\rm If $V=\bR$, we denote
$\bD_{\beta}^{1,p}(|\cH_{V}|)=\bD_{\beta}^{1,p}(|\cH|)$,
$\bL_{H,\beta}^{1,p}(V)=\bL_{H,\beta}^{1,p}$, and
$\widetilde{\bL}_{H,\beta}^{1,p}(V)=\widetilde{\bL}_{H,\beta}^{1,p}$
.}
\end{remark}

Note that the space $\bD_{\beta}^{1,p}(|\cH_{V}|$ is {\em not} the
particular instance of the space $\bD_{B}^{1,p}(|\cH_{V}|)$
(introduced in Section 2) obtained for $V=\bR$. The fundamental
difference between the two spaces is that
$\bD_{\beta}^{1,p}(|\cH_{V}|)$ contains $V$-valued random processes
$g=\{g(s,\cdot)\}_{s \in [0,T]}$, for an arbitrary Banach space $V$
(which has nothing to do with the underlying Hilbert space $\bR$ of
the fBm $\beta$), whereas the space $\bD_{B}^{1,p}(|\cH_{V}|$
contains $V$-valued random processes $U=\{U(s,\cdot)\}_{s \in
[0,T]}$, where $V$ {\em is} the underlying space of the Gaussian
process $B$.

In the present article, we let $V=H_p^n$.
Since $C_0^{\infty}$ is dense in $H_p^n$, we introduce the set
$\cS_{\beta}(\cE_{C_0^{\infty}})$ of smooth elementary processes of
the form
$$g(t,\cdot)=\sum_{i=1}^{m}F_i 1_{(t_{i-1},t_{i}]}(t)\phi_i(\cdot),
\quad t \in [0,T]$$ with $F_i \in \cS_{\beta}$, $0 \leq t_0<\ldots
<t_{m} \leq T$ and $\phi_i \in C_0^{\infty}$. The set
$\cS_{\beta}(\cE_{C_0^{\infty}})$ is dense in
$\bD_{\beta}^{1,p}(|\cH_{H_p^n}|)$ with respect to the norm $\|\cdot
\|_{\bD_{\beta}^{1,p}(|\cH_{H_p^n}|)}$. The space
$\widetilde{\bL}_{H,\beta}^{1,p}(H_p^n)$ is the completion of
$\cS_{\beta}(\cE_{C_0^{\infty}})$ in
$\bD_{\beta}^{1,p}(|\cH_{H_p^n}|)$, with respect to the norm
$\|\cdot\|_{{\bL}_{H,\beta}^{1,p}(H_p^n)}$. From
(\ref{ineq-norms-L-D-U}), it follows that
$\widetilde{\bL}_{H,\beta}^{1,p}(H_p^n) \subset
{\bL}_{H,\beta}^{1,p}(H_p^n)$.

For any $g \in \bL_{H,\beta}^{1,p}(H_p^n)$, we have:
\begin{equation}
\label{def-norm-in-LpH}
 \|g\|_{\bL_{H,\beta}^{1,p}(H_p^n)}^p=
\|g\|_{\bH_p^n}^p+\|D^{\beta}g\|_{\bH_{p,H}^n}^p,
\end{equation}
where
\begin{eqnarray*}
\bH_p^n&:=&L_p(\Omega \times [0,T],\cF \times \cB([0,T]);
H_p^n) \\
\bH_{p,H}^n&:=&L_p(\Omega \times [0,T], \cF \times \cB([0,T]);
L_{1/H}([0,T];H_p^n)).
\end{eqnarray*}

\vspace{2mm}

For an arbitrary element $g \in \bD^{1,p}_{\beta}(|\cH_{H_p^n}|)$,
we write $g(*,\cdot)=\{g(s,\cdot)\}_{s \in [0,T]}$.
Using (\ref{norm-of-g-Hpn}), for any $g \in
\bD_{\beta}^{1,p}(|\cH_{H_p^n}|)$ and $\phi \in C_0^{\infty}$, we
have:
\begin{eqnarray}
\label{first-estimate-H-V}
E\|(g(*,\cdot), \phi) \|_{|\cH|}^2 & \leq & N E\|g \|_{|\cH_{H_p^n}|}^2 \\
\label{second-estimate-H-V}
 E\|(D^{\beta}g(*,\cdot), \phi) \|_{|\cH|
\otimes |\cH|}^2 & \leq & N E\|D^{\beta}g \|_{|\cH| \otimes
|\cH_{H_p^n}|}^2,
\end{eqnarray}
where $N$ is a constant depending on $n,p$ and $\phi$.


\begin{proposition}
\label{fundamental-proposition}  a) If $g \in
\bD_{\beta}^{1,p}(|\cH_{H_p^n}|)$, then for any $\phi \in
C_0^{\infty}$, $(g(*,\cdot),\phi) \in \bD_{\beta}^{1,2}(|\cH|)$,
$D^{\beta}(g(*,\cdot),\phi)=(D^{\beta}g(*,\cdot),\phi)$, and
\begin{equation}
\label{Malliavin-estimate-H-V} \|(g(*,\cdot),\phi)
\|_{\bD_{\beta}^{1,p}(|\cH|)} \leq N \|g
\|_{\bD_{\beta}^{1,p}(|\cH_{H_p^n}|)},
\end{equation}
where $N$ is a constant depending on $n,p$ and $\phi$.

b) If $g\in \bL_{H,\beta}^{1,p}(H_p^n)$, then for any $\phi \in
C_0^{\infty}$, $(g(*,\cdot),\phi) \in \bL_{H,\beta}^{1,p}$, and
\begin{equation}
\label{ineq-norms-L-L(H)}
\|(g(*,\cdot),\phi)\|_{\bL_{H,\beta}^{1,p}} \leq N
\|g\|_{\bL_{H,\beta}^{1,p}(H_p^n)},
\end{equation} where
$N$ is a constant depending on $n,p$ and $\phi$.
\end{proposition}

\noindent {\bf Proof:} a) Using an approximation argument and the
completeness of the space $\bD_{\beta}^{1,p}(|\cH|)$, it suffices to
assume that $g(t,\cdot)=\sum_{i=1}^{m}F_i 1_{(t_i,t_{i+1}]}(t)
\phi_i$ with $F_i \in \cS_{\beta}$, $0 \leq t_1 < \ldots < t_{m+1}
\leq T$ and $\phi_i \in C_0^{\infty}$. Clearly,
$(g(*,\cdot),\phi)=\sum_{i=1}^{m} F_i (\phi_i,\phi)
1_{(t_i,t_{i+1}]} \in \cS_{\beta}(\cE) \subset
\bD_{\beta}^{1,2}(|\cH|)$, and due to the linearity of $D^{\beta}$,
$$D_t^{\beta}(g(s,\cdot),\phi)=\sum_{i=1}^{m}(D_t^{\beta}F_i)
(\phi_i,\phi)1_{(t_i,t_{i+1}]}(s)=(D_t^{\beta}g(s, \cdot),\phi).$$
Using (\ref{first-estimate-H-V}) and (\ref{second-estimate-H-V}), we
get:
\begin{eqnarray*}
\|(g(*,\cdot),\phi)\|_{\bD_{\beta}^{1,p}(|\cH|)}^p & =&
E\|(g(*,\cdot),\phi)\|_{|\cH|}^p+E\|D^{\beta}(g(*,\cdot),\phi)\|_{|\cH|
\otimes |\cH|}^p \\
& \leq & N (E\|g\|_{|\cH_{H_p^n}|}^p + E\|D^{\beta}g\|_{|\cH|
\otimes |\cH_{H_p^n}|}^p) = N \|g
\|_{\bD_{\beta}^{1,p}(|\cH_{H_p^n}|)}^{p}.
\end{eqnarray*}

b) By part a), $(g(*,\cdot),\phi) \in \bD_{\beta}^{1,p}(|\cH|)$.
Using (\ref{norm-of-g-Hpn}),
\begin{eqnarray*}
\|(g(*,\cdot),\phi)\|_{\bL_{H,\beta}^{1,p}}^p &=& E\int_0^T
|(g(s,\cdot),\phi)|^p ds +E \int_0^T \left( \int_0^T
|(D_t^{\beta}g(s,\cdot),\phi)|^{1/H}dt\right)^{pH}ds \\
& \leq & N E\int_0^T \|g(s,\cdot)\|_{H_p^n}^pds + E \int_0^T \left(
\int_0^T \|D_t^{\beta}g(s,\cdot)\|_{H_p^n}^{1/H}dt\right)^{pH}ds \\
&=& \|g\|_{\bL_{H,\beta}^{1,p}(H_p^n)}^p<\infty.
\end{eqnarray*}
$\Box$

\subsection{The case of a sequence of fBm's}
\label{sequence-fBms}

For any $p>1$ and $n \in \bR$, we let $H_p^n(l_2)$ be the set of all
sequences $u=(u^k)_k$ such that $u^k \in H_p^n$ for all $k$, and
$\|u\|_{H_p^n(l_2)}:=\| \ |(1-\Delta)^{n/2}u|_{l_2}
\|_{L_p}<\infty$.
By Minkowski's inequality, $\|u\|_{H_p^n(l_2)}^{2} \leq
\sum_{k=1}^{\infty} \|u^k\|_{H_p^n}^{2}$ (with equality if $p=2$).
By H\"older's inequality, for any $u \in H_p^n(l_2)$ and $\phi \in
C_0^{\infty}$, we have:
\begin{equation}
\label{ineq-norms-Hpn-l2} \sum_{k=1}^{\infty} |(u^k,\phi)|^2 \leq N
\|u\|_{H_p^n(l_2)}^{2}
\end{equation}
where $N$ is the same constant as in (\ref{norm-of-g-Hpn}).

Let $\beta^k=(\beta_t^k)_{t \in [0,T]}, k \geq 1$ be a sequence of
i.i.d. fBm's with Hurst index $H>1/2$, defined on the same
probability space $(\Omega, \cF, P)$. We first define the
$l_2$-analogue of the space $\bL_{H,\beta}^{1,p}$, introduced in
subsection \ref{single-fBm}.

\begin{definition}
\label{defin-of-L-1p} For any $p>1$, we denote by
$\bL_{H}^{1,p}(l_2)$ the set of all elements $u=(u^k)_k$ such that
$u^k \in \bD_{\beta^k}^{1,p}(|\cH|)$ for all $k$, and
$\|u\|_{\bL_{H}^{1,p}(l_2)}<\infty$, where
$$\|u\|_{{\bL}_{H}^{1,p}(l_2)}^{p}:=E \int_0^T
|u_s|_{l_2}^{p}ds+ E \int_0^T \left( \int_0^T
|D_{\theta}u_s|_{l_2}^{1/H} d\theta\right)^{pH}ds.$$

\end{definition}

The next lemma shows that condition (\ref{cond-on-uk}) in Theorem
\ref{connection} is satisfied for any $u=(u^k)_k \in
\bL_{H}^{1,p}(l_2)$.

\begin{lemma}
\label{remark-def-LH-l2} If  $p \geq 2$ and $u=(u^k)_k \in
\bL_{H}^{1,p}(l_2)$, then
$\sum_{k=1}^{\infty}\|u^k\|_{\bD_{\beta^k}^{1,2}(\cH)}^2<\infty$.
\end{lemma}

\noindent {\bf Proof:} Note that $\bD_{\beta^k}^{1,p}(|\cH|) \subset
\bD_{\beta^k}^{1,2}(|\cH|)$. For any $u \in \bL_{H}^{1,p}(l_2)$, we
have:
\begin{eqnarray*}
\sum_{k=1}^{\infty}\|u^k\|_{\bD_{\beta^k}^{1,2}(\cH)}^2 & \leq &
\sum_{k=1}^{\infty} \left\{E\int_0^T |u_s^k|^2 ds+E \int_0^T \left(
\int_0^T |D_{\theta}^{\beta^k}u_s^k|^{1/H}d \theta
\right)^{2H}ds\right\} \\
& \leq & E \int_0^T |u_s|_{l_2}^2 ds+ E \int_0^T \left(\int_0^T
|D_{\theta} u_s|_{l_2}^{1/H} d \theta \right)^{2H}ds \\
& \leq & C_{p,H,T}\|u\|_{\bL_{H}^{1,p}(l_2)}^p<\infty,
\end{eqnarray*}
where $C_{p,H,T}$ is a constant depending on $p,H$ and $T$. The
first inequality above is due to (\ref{ineq-norms-1}) and
(\ref{ineq-norms-2}), the second is due to Minkowski's inequality,
and the third is due to H\"{o}lder's inequality. $\Box$

We now introduce the definition of the space
$\widetilde{\bL}_{H}^{1,p}(H_p^n,l_2)$, in which we are allowed to
select the coefficients $(g^k)_k$ multiplying the noise in the
stochastic heat equation.

\begin{definition} Let $p >1$ be arbitrary.

a) We denote by $\bL_{H}^{1,p}(H_p^n,l_2)$ the set of all elements
$g=(g^k)_k$ such that $g^k \in \bD_{\beta^k}^{1,p}(|\cH_{H_p^n}|)$
for all $k$, and $\|g\|_{\bL_{H}^{1,p}(H_p^n,l_2)}<\infty$, where
$$\|g\|_{\bL_{H}^{1,p}(H_p^n,l_2)}^{p}:=E \int_0^T
|g(s, \cdot)|_{H_p^n(l_2)}^{p}ds+ E \int_0^T \left( \int_0^T
|D_{\theta}g(s, \cdot)|_{H_p^n(l_2)}^{1/H}d\theta\right)^{pH}ds.$$

b) We let $\widetilde{\bL}_{H}^{1,p}(H_p^n,l_2)$ be the set of all
$g \in \bL_{H}^{1,p}(H_p^n,l_2)$ for which there exists a sequence
$(g_j)_j \subset \bL_{H}^{1,p}(H_p^n,l_2)$ such that
$\|g_j-g\|_{\bL_{H}^{1,p}(H_p^n,l_2)} \to 0$ as $j \to \infty$,
$g_j^k=0$ for $k>K_j$, and $g_j^k \in
\cS_{\beta^k}(\cE_{C_0^{\infty}})$ for $k \leq K_j$, i.e.
$$g_j^k(t,\cdot)=\sum_{i=1}^{m_{jk}} F_{i}^{jk}
1_{(t_{i-1}^{jk},t_{i}^{jk}]}(t) \phi_{i}^{jk}(\cdot), \quad t \in
[0,T],$$ with $F_{i}^{jk} \in \cS_{\beta^k}$, $0 \leq t_0^{jk}<
\ldots < t_{m_{jk}}^{jk} \leq T$ (non-random) and $\phi_i^{jk} \in
C_0^{\infty}$.
\end{definition}

Note that, for any $g \in \bL_{H}^{1,p}(H_p^n,l_2)$,
\begin{equation}
\label{def-norm-in-LpH-l2}
 \|g\|_{\bL_{H}^{1,p}(H_p^n,l_2)}^p=
\|g\|_{\bH_p^n(l_2)}^p+\|D g\|_{\bH_{p,H}^n(l_2)}^p,
\end{equation}
where
\begin{eqnarray*}
\bH_p^n(l_2)&:=&L_p(\Omega \times [0,T],\cF \times \cB([0,T]);
H_p^n(l_2)) \\
\bH_{p,H}^n(l_2)&:=&L_p(\Omega \times [0,T], \cF \times \cB([0,T]);
L_{1/H}([0,T];H_p^n(l_2))).
\end{eqnarray*}

\begin{lemma}
\label{g-gk-in-L} If  $g=(g^k)_k \in \bL_{H}^{1,p}(H_p^n,l_2)$, then
$g^k \in \bL_{H,\beta^k}^{1,p}(H_p^n)$ for all $k$, and
$$\|g^k\|_{\bL_{H,\beta^k}^{1,p}(H_p^n)} \leq
\|g\|_{\bL_{H}^{1,p}(H_p^n,l_2)} \quad \mbox{for all} \ k.$$ In
particular, if $g=(g^k)_k \in \widetilde{\bL}_{H}^{1,p}(H_p^n,l_2)$,
then $g^k \in \widetilde{\bL}_{H,\beta^k}^{1,p}(H_p^n)$ for all $k$.
\end{lemma}

\noindent {\bf Proof:} We have:
\begin{eqnarray*}
\lefteqn{ \|g^k\|_{\bL_{H,\beta^k}^{1,p}(H_p^n)}^{p} = E \int_0^T
\|g^k(s, \cdot)\|_{H_p^n}^p ds+ E \int_0^T \left(\int_0^T
\|D_{\theta}^{\beta^k}g^k(s, \cdot) \|_{H_p^n}^{1/H}\right)^{pH}ds }
\\
& & = E \int_0^T \|(1-\Delta)^{n/2}g^k(s, \cdot)\|_{L_p}^p ds+ E
\int_0^T \left(
\int_0^T\|D_{\theta}^{\beta^k}[(1-\Delta)^{n/2}g^k(s, \cdot)] \
\|_{H_p^n}^{1/H}\right)^{pH}ds \\
& & \leq E \int_0^T \|\ |(1-\Delta)^{n/2}g(s,
\cdot)|_{l_2}\|_{L_p}^p ds+ E \int_0^T \left( \int_0^T\| \ |
D_{\theta}[(1-\Delta)^{n/2}g(s, \cdot)] \ |_{l_2}
\|_{H_p^n}^{1/H}\right)^{pH}ds
\\
& & = \|g \|_{\bL_{H}^{1,p}(H_p^n,l_2)}^{p}.
\end{eqnarray*}
The second statement follows from the definitions of spaces
$\widetilde{\bL}_{H}^{1,p}(H_p^n,l_2)$ and
$\widetilde{\bL}_{H,\beta^k}^{1,p}(H_p^n)$. $\Box$

\section{The Main Result}

The following definition introduces the solution space (see
Definition 3.1 of \cite{krylov99}).

\begin{definition}
\label{def-space-cH} Let $p \geq 2$ be arbitrary.

Let $u=\{u(t,\cdot)\}_{t \in [0,T]}$ be a $\cD$-valued random
process defined on the probability space $(\Omega, \cF,P)$. We write
$u \in \cH_{p,H}^n$ if:

(i) $u(0,\cdot) \in L_p(\Omega,\cF,H_p^{n-2/p})$;

(ii) $u \in \bH_p^{n}$, $u_{xx} \in \bH_p^{n-2}$;

(iii) there exist $f \in \bH_p^{n-2}$ and $g \in
\widetilde{\bL}_{H}^{1,p}(H_p^{n-1},l_2)$ such that for any $\phi
\in C_0^{\infty}$, the equality
\begin{equation}
\label{def-sol-space} (u(t,\cdot),\phi)=(u(0,\cdot),\phi)+\int_0^t
(f(s,\cdot),\phi)ds+\sum_{k=1}^{\infty}\int_0^t
(g^k(s,\cdot),\phi)\delta \beta_s^k
\end{equation}
holds for any $t \in [0,T]$ a.s. We define
\begin{equation}
\label{norm-in-cH}
\|u\|_{\cH_{p,H}^n}=(E\|u(0,\cdot)\|_{H_p^{n-2/p}}^p)^{1/p}+
\|u_{xx}\|_{\bH_{p}^{n-2}}+
\|f\|_{\bH_p^{n-2}}+\|g\|_{\bL_{H}^{1,p}(H_p^{n-1},l_2)}.
\end{equation}
If $u \in \cH_{p,H}^n$, we write ${\bf D}u:=f$, ${\bf S} u:=g$ and
$du=f dt+\sum_{k=1}^{\infty} g^k \delta \beta^k_t$, $t \in [0,T]$.

We say that $u \in \cH_{p,H}^{n}$ is a {\em solution} of
(\ref{heat}) if ${\bf D}u=\Delta u+f$ and ${\bf S}u=g$.
\end{definition}

\begin{remark}
{\rm The series of stochastic integrals in (\ref{def-sol-space})
converges uniformly in $t$, in probability. More precisely, if $g
\in \bL_{H}^{1,p}(H_p^{n},l_2), \phi \in C_0^{\infty}$ are
arbitrary, and we let $u_t^k=(g^k(t,\cdot),\phi),t \in [0,T]$, then
$$u \in {\bL}_{H}^{1,p}(l_2).$$ (To see this, note that by Lemma
\ref{g-gk-in-L}, $g^k \in \bL_{H,\beta^k}^{1,p}(H_p^n)$ for all $k$.
By Proposition \ref{fundamental-proposition}, $u^k \in
\bL_{H,\beta^k}^{1,p}$ for all $k$. Since by
(\ref{ineq-norms-Hpn-l2}), $|u_s|_{l_2} \leq N \|g(s,
\cdot)\|_{H_p^n(l_2)}$ and $|D_{\theta}u_s|_{l_2} \leq N
\|D_{\theta}g(s, \cdot)\|_{H_p^n(l_2)}$, we get:
$\|u\|_{{\bL}_{H}^{1,p}(l_2)} \leq N
\|g\|_{{\bL}_{H}^{1,p}(H_p^n,l_2)}< \infty$.) By Lemma
\ref{remark-def-LH-l2}, $\sum_{k=1}^{\infty} \|u^k
\|_{\bD_{\beta^k}^{1,2}(|\cH|)}^2<\infty$. Denoting
$X_t^{(N)}:=\sum_{k=1}^{N}\int_{0}^{t} u_s^k\delta \beta_s^k$ and
$X_{t}:=\sum_{k=1}^{\infty}\int_{0}^{t} u_s^k\delta \beta_s^k$,
relation (\ref{convergence-in-sup-metric}) shows that
$$\lim_{N\to \infty}P(\sup_{t \leq T}|X_t^{(N)}-X_t|\geq \varepsilon)=0, \quad
\mbox{for any} \ \varepsilon>0.$$

\noindent In what follows, we work with an a.s. continuous
modification of $X=(X_t)_{t \in [0,T]}$.}
\end{remark}

\begin{remark}
{\rm By the definition of the norm in $\cH_{p,H}^n$, the operators
${\bf D}: \cH_{p,H}^n \to \bH_{p}^{n-2}([0,T])$ and ${\bf S}:
\cH_{p,H}^n \to \widetilde{\bL}_{H}^{1,p}(H_p^{n-1},l_2)$ are
continuous.}
\end{remark}

\begin{proposition}
\label{(1-Delta)-isometry} (a) The operator $(1-\Delta)^{m/2}$ maps
isometrically $\widetilde{\bL}_{H}^{1,p}(H_p^n,l_2)$  onto
$\widetilde{\bL}_{H}^{1,p}(H_p^{n-m},l_2)$.

(b) The operator $(1-\Delta)^{m/2}$ maps isometrically
$\cH_{p,H}^{n}$ onto  $\cH_{p,H}^{n-m}$.
\end{proposition}

\noindent {\bf Proof:} (a) By the definition of
$\widetilde{\bL}_{H}^{1,p}(H_p^n,l_2)$, it suffices to prove that
$(1-\Delta)^{m/2}$ maps isometrically
$\widetilde{\bL}_{H,\beta}^{1,p}(H_p^n)$  onto
$\widetilde{\bL}_{H,\beta}^{1,p}(H_p^{n-m})$, for a fixed fBm
$\beta=(\beta_t)_{t \in [0,T]}$.

Let $g \in \bL_{H,\beta}^{1,p}(H_p^n)$ be arbitrary.
By Proposition \ref{fundamental-proposition},
\begin{eqnarray*}
(D^{\beta}[(1-\Delta)^{m/2}g(*,\cdot)],\phi) & =&
D^{\beta}((1-\Delta)^{m/2}g(*,\cdot),\phi)=
D^{\beta}(g(*,\cdot),(1-\Delta)^{m/2}\phi) \\
& = & (D^{\beta}g(*,\cdot),(1-\Delta)^{m/2}\phi)=
((1-\Delta)^{m/2}[D^{\beta}g(*,\cdot)],\phi),
\end{eqnarray*}
for any $\phi \in C_0^{\infty}$, i.e.
$$
D_t^{\beta}[(1-\Delta)^{m/2}g(s,\cdot)]=
(1-\Delta)^{m/2}[D_t^{\beta}g(s,\cdot)], \quad \forall s,t \in
[0,T].
$$

\noindent Using an approximation argument and the fact that
$\|u\|_{H_p^{n}}=\|(1-\Delta)^{m/2}u\|_{H_p^{n-m}}$ for any $u \in
H_p^n$, we conclude that $(1-\Delta)^{m/2}g \in
\bD_{\beta}^{1,p}(|\cH_{H_p^{n-m}}|)$ and
\begin{eqnarray*}
 \|(1-\Delta)^{m/2}g\|_{\bL_{H,\beta}^{1,p}(H_p^{n-m})}^p&=&
\|(1-\Delta)^{m/2}g\|_{\bH_p^{n-m}}^p+
\|D^{\beta}[(1-\Delta)^{m/2}g] \|_{\bH_{p,H}^{n-m}}^p \\
&=&\|g\|_{\bH_p^n}^p+ \|D^{\beta}g\|_{\bH_{p,H}^n}^p =
\|g\|_{\bL_{H,\beta}^{1,p}(H_p^{n})}^p<\infty.
\end{eqnarray*}
This proves that $(1-\Delta)^{m/2}g \in
\bL_{H,\beta}^{1,p}(H_p^{n-m})$. Finally, if $g \in
\widetilde{\bL}_{H,\beta}^{1,p}(H_p^n)$, then an approximation
argument shows that $(1-\Delta)^{m/2}g \in
\widetilde{\bL}_{H,\beta}^{1,p}(H_p^{n-m})$.


(b) This is a consequence of part a). See Remark 3.8 of
\cite{krylov99}. $\Box$

\vspace{2mm}

\begin{theorem}
\label{krylov-theorem3-7} (a) If $u \in \cH_{p,H}^{n}$, then $u \in
C([0,T],H_p^{n-2})$ a.s.,
$$E \sup_{t \leq T} \|u(t,\cdot)\|_{H_p^{n-2}}^{p} \leq N
 \|u\|_{\cH_{p,H}^{n}}^{p} \quad
and \quad \|u\|_{\bH_p^{n}} \leq N
 \|u\|_{\cH_{p,H}^{n}},$$
where $N$ is a constant which depends on $p,H,T$ and $d$.

\noindent (b) $\cH_{p,H}^{n}$ is a Banach space with the norm
(\ref{norm-in-cH}).
\end{theorem}

\noindent {\bf Proof:} (a) By Proposition \ref{(1-Delta)-isometry},
it suffices to take $n=0$. We use the same argument as in the proof
of Theorem 3.7 of \cite{krylov99}. We refer the reader to this proof
for the notation. In our case, we only need to justify that:
$$E \sup_{t \leq
T}\left\|\sum_{k=1}^{\infty}\int_0^t
g^{(\varepsilon)k}(s,\cdot)\delta \beta_s^k \right\|_{L_p}^p \leq C
\|u\|_{\cH_{p,H}^{2}}^{p},$$ where $C$ is a constant which depends
on $p,H$ and $T$.

Using Corollary \ref{max-inequality-p-cor}, for any $x \in \bR^d$,
we have:
\begin{eqnarray*}
 E\sup_{t \leq T} \left|\sum_{k=1}^{\infty} \int_0^t
g^{(\varepsilon)k}(s,x) \delta \beta_s^k \right|^p & \leq & C
\left\{ E\int_0^T |g^{(\varepsilon)}(s,x)|_{l_2}^{p}ds
 + \right.  \\
 & & \left. E
\int_0^T \left( \int_0^T |D_{\theta}
g^{(\varepsilon)}(s,x)|_{l_2}^{1/H} d\theta \right)^{pH}ds
 \right\},
\end{eqnarray*}

\noindent where $C$ is a constant depending on $p,H$ and $T$. We
integrate with respect to $x$. Using Minkowski's inequality
 and the fact that
$\|h^{(\varepsilon)}\|_{L_2} \leq \|h\|_{L_2}$ for any $h \in L_2$,
we get:
\begin{eqnarray*}
\lefteqn{E \sup_{t \leq T}\left\|\sum_{k=1}^{\infty}\int_0^t
g^{(\varepsilon)k}(s,\cdot)\delta \beta_s^k \right\|_{L_p}^p  \leq C
\left\{E \int_0^T  \int_{\bR^d} |g^{(\varepsilon)}(s,x)|_{l_2}^{p}dx
ds+ \right. } \\
& & \left. E \int_0^T \int_{\bR^d} \left(\int_0^T
|D_{\theta}^{\beta^k}g^{(\varepsilon)}(s,x)|_{l_2}^{1/H}d\theta
\right)^{pH}dx ds \right\} \\
& & \leq C \left\{ E \int_0^T   \| \
|g^{(\varepsilon)}(s,\cdot)|_{l_2} \|_{L_p}^{p} ds + E \int_0^T
\left(\int_0^T \| \ | D_{\theta}g^{(\varepsilon)}(s, \cdot)|_{l_2}
\|_{L_p}^{1/H} d\theta \right)^{pH}ds \right\} \\
& & \leq C \left\{ E \int_0^T   \| \ |g(s,\cdot)|_{l_2} \|_{L_p}^{p}
ds + E \int_0^T \left(\int_0^T \| \ | D_{\theta}g(s, \cdot)|_{l_2}
\|_{L_p}^{1/H} d\theta \right)^{pH}ds \right\} \\
& & = C \left\{ E \int_0^T \|g(s, \cdot)\|_{L_p(l_2)}^{p}ds+ E
\int_0^T \left( \int_0^T \|D_{\theta} g(s, \cdot)
\|_{L_p(l_2)}^{1/H}d\theta \right)^{pH}ds  \right\}\\
& & = C \ \|g\|_{\bL_{H}^{1,p}(L_p,l_2)}^{p} \leq  C
\|g\|_{\bL_{H}^{1,p}(H_p^1,l_2)}^{p} \leq C
\|u\|_{\cH_{p,H}^{2}}^{p}.
\end{eqnarray*}

(b) Let $\{u_j\}_{j}$ be a Cauchy sequence in $\cH_{p,H}^{n}$. By
(a), $\{u_j\}_{j}$ is a Cauchy sequence in $\bH_{p}^{n}$. Hence,
there exists $u \in \bH_{p}^{n}$ such that $\|u_j-u\|_{\bH_{p}^{n}}
\to 0$. Moreover, $u_{xx} \in \bH_{p}^{n-2}$ and
$\|u_{jxx}-u_{xx}\|_{\bH_{p}^{n-2}} \to 0$.

Say $u_j$ satisfies (\ref{def-sol-space}) for $f_j \in \bH_p^{n-2}$,
$g_j \in \widetilde{\bL}_{H}^{1,p}(H_p^{n-1},l_2)$: for any $\phi
\in C_0^{\infty}$,
\begin{equation}
\label{def-sol-space-j}
(u_j(t,\cdot),\phi)=(u_j(0,\cdot),\phi)+\int_0^t
(f_j(s,\cdot),\phi)ds+\sum_{k=1}^{\infty}\int_0^t
(g_j^k(s,\cdot),\phi)\delta \beta_s^k
\end{equation}
 for any $t \in [0,T]$ a.s.
 Then $\{u_j(0,\cdot)\}_j,
 \{f_j\}_j$ and $\{g_j\}_j$ are Cauchy in the
 (complete) spaces $L_p(\Omega,\cF;H_p^{n-2/p})$, $\bH_p^{n-2}$ and
 $\widetilde{\bL}_{H}^{1,p}(H_p^{n-1},l_2)$, respectively. Hence,
 there exist
 $u(0,\cdot) \in L_p(\Omega,\cF,H_p^{n-2/p}),f \in
\bH_p^{n-2}$, $g \in \widetilde{\bL}_{H}^{1,p}(H_p^{n-1},l_2)$ such
that $E\|u_j(0,\cdot)-u(0,\cdot)\|_{H_p^{n-2/p}} \to 0$,
$\|f_j-f\|_{\bH_p^{n-2}} \to 0$ and
$\|g_j-g\|_{\bL_{H}^{1,p}(H_p^{n-1},l_2)} \to 0$.

Since $\|u_j-u\|_{\bH_{p}^{n}} \to 0$, there exists a subsequence of
indices $j$ such that $\|u_j(t,\cdot)-u(t,\cdot)\|_{H_p^{n}} \to 0$
a.e. in $(\omega,t)$. Say that this happens for $\omega \in \Omega
\verb2\2 \Gamma$ and $t \in [0,T] \verb2\2 U$, where $\Gamma,U$ are
negligible sets.

 Fix $t \in [0,T] \verb2\2 U$. We are now passing to the limit in (\ref{def-sol-space-j}). On the
left hand side, $|(u_j(t,\cdot)-u(t,\cdot),\phi)| \leq N
\|u_j(t,\cdot)-u(t,\cdot) \|_{H_p^n} \to 0$ a.s. On the right hand
side of (\ref{def-sol-space-j}), the first two terms clearly
converge to $(u(0,\cdot),\phi)$, respectively $\int_0^t
(f(s,\cdot),\phi)ds$. For the third term, we invoke Corollary
\ref{max-inequality-p-cor} and (\ref{ineq-norms-L-L(H)}):
\begin{eqnarray*}
E\left|\sum_{k=1}^{\infty}\int_0^t
(g_j^k(s,\cdot)-g^k(s,\cdot),\phi) \delta \beta_s^k \right|^p & \leq
& N \|(g_j^k(*,\cdot)-g^k(*,\cdot),\phi)\|_{\bL_{H}^{1,p}(l_2)}^p
\\
& \leq & N\|g_j^k-g^k\|_{\bL_{H}^{1,p}(H_p^{n-1},l_2)}^p    \to 0, \
{\rm as} \ j \to \infty.
\end{eqnarray*}
Therefore, $\sum_{k=1}^{\infty}\int_0^t
(g_j^k(s,\cdot)-g^k(s,\cdot),\phi) \delta \beta_s^k  \to 0$ a.s.
(for a subsequence of indices $j$). We infer that for every $\phi
\in C_0^{\infty}$ and for any $t \in [0,T] \verb2\2 U$, equality
(\ref{def-sol-space}) holds almost surely (with the negligible set
depending on $t$).

To conclude that $u \in \cH_{p,H}^{n}$, it remains to show that
equality (\ref{def-sol-space}) holds for any $t \leq T$ a.s. (i.e.
the negligible set does not depend on $t$). For this, it suffices to
note that the process $(u(*,\cdot),\phi)$ is continuous a.s. This
follows from the a.s. continuity of processes $(u_j(*,\cdot),\phi)$,
by noting that $(u_j(t,\cdot),\phi)$ converges to
$(u(t,\cdot),\phi)$ uniformly in $t$, in probability.
$\Box$

\vspace{3mm}

The next theorem is the main result of the present article.

\begin{theorem}
\label{main-theorem} Let $p \geq 2$ and $n \in \bR$ be arbitrary.
Let
$$f \in \bH_p^{n-2}, \quad g \in
\widetilde{\bL}_{H}^{1,p}(H_p^{n-1},l_2) \quad \mbox{and} \quad u_0
\in L_p(\Omega, \cF, H_p^{n-2/p}).$$ Then the Cauchy problem for
equation (\ref{heat}) with initial condition $u(0,\cdot)=u_0$ has a
unique solution $u \in \cH^{n}_{p,H}$. For this solution, we have
\begin{equation}
\label{krylov-relation4-25} \|u\|_{\cH_{p,H}^{n}} \leq N
\{\|f\|_{\bH_{p}^{n-2}([0,T])}+
\|g\|_{\bL_{H}^{1,p}(H_p^{n-1},l_2)}+(E\|u_0\|_{H_p^{n-2/p}}^{p})^{1/p}\},
\end{equation}
where $N$ is a constant depending on $p,d,T$ and $H$.
\end{theorem}

\noindent {\bf Proof:} We first prove that it suffices to take
$u_0=0$. To see this, we assume without loss of generality that
$n=2$ (using Proposition \ref{(1-Delta)-isometry}). By Theorem 2.1
of \cite{krylov99}, for every $\omega \in \Omega$ fixed, the
equation $du=\Delta u \ dt$ with initial condition $u_0$ has a
unique solution $\bar{u} \in H_{p}^{1,2}$, and
$\|\bar{u}\|_{H_p^{1,2}} \leq N \|u_0\|_{H_p^{2-2/p}}$ and
$\|\bar{u}_{xx}\|_{L_p((0,T) \times \bR^d)}
\leq N \|u_0\|_{H_p^{2-2/p}}$. From here, one can show that $\bar{u}
\in \cH_{p,H}^2$ and $\|\bar{u}\|_{\cH_{p,H}^2} \leq N
\|u_0\|_{H_p^{2-2/p}}$. Suppose that equation (\ref{heat}) with zero
initial condition has a unique solution $v \in \cH_{p,H}^2$, and
$\|v\|_{\cH_{p,H}^2} \leq N
(\|f\|_{\bH_{p}^0}+\|g\|_{\bL_{H}^{1,p}(H_p^1,l_2)})$. Then
$u:=v+\bar{u} \in \cH_{p,H}^2$ is a solution of (\ref{heat}) with
initial condition $u_0$, and (\ref{krylov-relation4-25}) holds.

For the remaining part of the proof, we assume that $u_0=0$. By
Proposition \ref{(1-Delta)-isometry}, it is enough to consider only
one particular value of $n$. We take $n=1$.

{\em Case 1.} Suppose that $g^k=0$ for $k>K$, and
$$g^k(t,\cdot)=
\sum_{i=1}^{m_k}F_i^{k}1_{(t_{i-1}^{k},t_{i}^{k}]}(t)g_i^{k}(\cdot),
\quad t \in [0,T], \ k \leq K,$$ where $F_{i}^k \in \cS_{\beta^k}$,
$0 \leq t_{0}^k < \ldots < t_{m_k}^k \leq T$, and $g_i^{k} \in
C_0^{\infty}$.

Let $v(t,x)=\sum_{k=1}^{\infty} \int_0^t g^k(s,x) \delta \beta_s^k$
and $z(t,x)=\int_0^t T_{t-s}(\Delta v+f)(s,\cdot)(x)ds$. One can
show that $u=v+z$ is a solution of (\ref{heat}).

Let $u_1(t,x)=\int_0^t T_{t-s}[f(s,\cdot)](x)ds$. We first show that
\begin{equation}
\label{norm-uu1} \|u-u_1 \|_{\bH_p^0([0,T])} \leq  N
\|g\|_{\bL_{H}^{1,p}(L_p,l_2)}, \quad
 \|u_{x}-u_{1x}
\|_{\bH_p^0([0,T])} \leq   N \|g\|_{\bL_{H}^{1,p}(L_p,l_2)},
\end{equation}
where $N$ is a constant depending on $p,d,T$ and $H$.

By definition, $ u(t,x)-u_1(t,x)=v(t,x)+\int_0^t T_{t-s} (\Delta
v)(s,\cdot)(x)ds$. Note that $v(s,x)= \sum_{k=1}^{\infty}
\sum_{i=1}^{m_k} g_i^k(x) \int_0^s F_i^k 1_{(t_{i-1}^{k}, t_i^k]}(r)
\delta \beta_r^k$.
Using the stochastic Fubini's theorem and the fact that $\int_r^t
T_{t-s}(\Delta g_i^k)(x)ds=T_{t-r}g_i^k(x)-g_i^k(x)$, we get:
\begin{equation}
\label{formula-for-u2} u(t,x)-u_1(t,x)=
\sum_{k=1}^{\infty}\sum_{i=1}^{m_k}\int_0^t F_i^k
1_{(t_{i-1}^{k},t_i^k]}(r)T_{t-r}g_i^k(x)\delta \beta_r^k=
\sum_{k=1}^{\infty}\int_0^t T_{t-r}g^{k}(r,\cdot)(x) \delta
\beta_r^k.
\end{equation}

By Corollary \ref{max-inequality-p-cor},
\begin{eqnarray}
\nonumber \lefteqn{\|u-u_1 \|_{\bH_p^0}^p = \int_0^T \int_{\bR^d} E
\left|\sum_{k=1}^{\infty}\int_0^t
T_{t-s}g^k(s,\cdot)(x)\delta \beta_s^k\right|^p dx dt } \\
\nonumber & & \leq C \left\{ \int_0^T \int_{\bR^d} E \int_0^t
\left(\sum_{k=1}^{\infty}|T_{t-s}g^k(s,\cdot)(x)|^{2} \right)^{p/2}
ds
 dx dt + \right. \\
 \nonumber
& & \left. \int_0^T  \int_{\bR^d} E \int_0^t   \left[ \int_0^T
\left(\sum_{k=1}^{\infty}
|D_{\theta}^{\beta^k}[T_{t-s}g^k(s,\cdot)(x)]|^2\right)^{1/(2H)}
d\theta
\right]^{pH} ds dx dt \right\} \\
\label{u-u1} & & :=C(I_1+I_2).
\end{eqnarray}

By Theorem \ref{max-inequality-p},
\begin{eqnarray}
\nonumber \lefteqn{\|u_{x}-u_{1x} \|_{\bH_p^0}^p = \int_0^T
\int_{\bR^d} E \left|\sum_{k=1}^{\infty}\int_0^t
T_{t-s}g_x^k(s,\cdot)(x)\delta \beta_s^k\right|^p dx dt } \\
\nonumber & & \leq C \left\{ \int_0^T \int_{\bR^d} E
\left( \int_0^t  \sum_{k=1}^{\infty}|T_{t-s}g_x^k(s,\cdot)(x)|^2 ds \right)^{p/2} dx dt + \right. \\
\nonumber & & \left. \int_0^T  \int_{\bR^d} E \left\{\int_0^t \left[
\int_0^T \left(\sum_{k=1}^{\infty}
|D_{\theta}^{\beta^k}[T_{t-s}g_x^k(s,\cdot)(x)]|^2\right)^{1/(2H)}
d\theta
\right]^{2H} ds \right\}^{p/2}dx dt \right\} \\
\label{ux-u1x} & & :=C(J_1+J_2).
\end{eqnarray}

For evaluating the terms $I_2$ and $J_2$ above, we need to observe
that:
\begin{equation}
\label{D-commutes-with-T} D_{\theta}^{\beta^k}
[T_{t-s}g^k(s,\cdot)(x)]=T_{t-s}[D_{\theta}^{\beta^k}g^k(s,\cdot)](x).
\end{equation}
(This is a consequence of Proposition
\ref{fundamental-proposition}.(a), and the fact that $T_{t-s}g^k(s,
\cdot)(x)=(g^k(s, \cdot)*G_{t-s})(x)=(g^k(s,\cdot),
G_{t-s}(x-\cdot))$.)

By (\ref{ineq-Tt-g-norm-p-V}) (see Appendix A) and Minkowski's
inequality, we have:
\begin{eqnarray}
\nonumber I_1 &=& E \int_0^T \int_0^t \int_{\bR^d} |T_{t-s}g(s,
\cdot)(x)|_{l_2}^{p}dx ds dt=E \int_0^T \int_0^t \|T_{t-s}g(s,
\cdot)\|_{L_p(l_2)}^{p}ds dt \\
\label{u-u1-term1} & \leq & E \int_0^T \int_0^t \|g(s,
\cdot)\|_{L_p(l_2)}^{p}ds dt \leq T \|g
\|_{\bH_{p}^{0}(l_2)}^{p} \\
\nonumber I_2 &=& E \int_0^T  \int_0^t \int_{\bR^d} \left( \int_0^T
|T_{t-s} [D_{\theta}g(s,\cdot)](x)|_{l_2}^{1/H} d\theta \right)^{pH}
dx ds dt
\\
\nonumber & \leq & E \int_0^T  \int_0^t \left[ \int_0^T
\left(\int_{\bR^d} |T_{t-s} [D_{\theta}g(s,\cdot)](x)|_{l_2}^{p} dx
\right)^{1/(pH)}
d\theta \right]^{pH} ds dt \\
\nonumber & = & E \int_0^T  \int_0^t \left( \int_0^T  \|T_{t-s}
[D_{\theta}g(s,\cdot)] \ \|_{L_p(l_2)}^{1/H}  d\theta
\right)^{pH} ds dt \\
\nonumber & \leq & E \int_0^T  \int_0^t \left( \int_0^T  \|
D_{\theta}g(s,\cdot)  \|_{L_p(l_2)}^{1/H}  d\theta
\right)^{pH} ds dt \\
\label{u-u1-term2} & \leq & T E \int_0^T  \left( \int_0^T  \|
D_{\theta}g(s,\cdot) \|_{L_p(l_2)}^{1/H}  d\theta \right)^{pH} ds =T
\|D g \|_{\bH_{p,H}^{0}(l_2)}^{p}.
\end{eqnarray}

\noindent From (\ref{u-u1}), (\ref{u-u1-term1}) and
(\ref{u-u1-term2}), we conclude that:
$$\|u-u_1 \|_{\bH_p^0}^p  \leq C T( \|g \|_{\bH_{p}^{0}(l_2)}^{p}
+ \|D g \|_{\bH_{p,H}^{0}(l_2)}^{p})=CT \|g
\|_{\bL_{H}^{1,p}(L_p,l_2)}^{p}.$$

Using Theorem \ref{theorem1-1-krylov} (Appendix A) and Minkowski's
inequality, we have:
\begin{eqnarray}
\nonumber J_1 &=& E \int_{\bR^d}\int_0^T \left( \int_0^t
|T_{t-s}g_x(s, \cdot)(x)|_{l_2}^{2} ds \right)^{p/2} dt dx \\
\label{ux-u1x-term1} & \leq & N E \int_{\bR^d} \int_0^T
|g(s,x)|_{l_2}^{p}ds dx=N \|g\|_{\bH_{p}^{0}(l_2)}^{p}
\end{eqnarray}

Using Theorem \ref{gener-Th1-1-Krylov} (Appendix A), we have:
\begin{eqnarray}
\nonumber J_2 &=& E  \int_{\bR^d}  \int_0^T  \left[\int_0^t \left(
\int_0^T |T_{t-s}[D_{\theta}g_x(s,\cdot)](x)|_{l_2}^{1/H} d\theta
\right)^{2H} ds \right]^{p/2}dt dx \\
 \nonumber & \leq & N E
\int_0^T \left[\int_0^T \left( \int_{\bR^d}|D_{\theta}
g(s,x)|_{l_2}^{p} dx \right)^{1/(pH)} d
\theta \right]^{pH} ds \\
\label{ux-u1x-term2}  &=& N E \int_0^T   \left(\int_0^T \|D_{\theta}
g(s,x)\|_{L_p(l_2)}^{1/H} d \theta \right)^{pH} ds= N \|D
g\|_{\bH_{p,H}^{0}(l_2)}^{p}.
\end{eqnarray}

\noindent From (\ref{ux-u1x}), (\ref{ux-u1x-term1}) and
(\ref{ux-u1x-term2}), we infer that:
$$\|u_{x}-u_{1x} \|_{\bH_p^0}^p \leq C N (\|g\|_{\bH_{p}^{0}(l_2)}^{p}+
\|D g\|_{\bH_{p,H}^{0}(l_2)}^{p})=C
N\|g\|_{\bL_{H}^{1,p}(L_p,l_2)}^{p}.$$ This concludes the proof of
(\ref{norm-uu1}).

It remains to prove that $u \in \cH_{p,H}^{1}$. Using
(\ref{norm-uu1}), we have:
\begin{eqnarray}
\label{prelim-est1}  \|u\|_{\bH_p^0} &\leq &  \|u_1\|_{\bH_p^0}+
\|u-u_1\|_{\bH_p^0}  \leq  N (
\|f\|_{\bH_{p}^{-1}}+\|g\|_{\bL_{H}^{1,p}(L_p,l_2)}) \\
 \nonumber
\|u_{xx}\|_{\bH_p^{-1}} &\leq & \|u_{1xx}\|_{\bH_p^{-1}}+
\|u_{xx}-u_{1xx}\|_{\bH_p^{-1}} \leq \|u_{1x}\|_{\bH_p^{0}}+
\|u_{x}-u_{1x}\|_{\bH_p^{0}}\\
\label{prelim-est2} & \leq & N (
\|f\|_{\bH_{p}^{-1}}+\|g\|_{\bL_{H}^{1,p}(L_p,l_2)}).
\end{eqnarray}

\noindent Using the fact that $\|\phi\|_{H_p^1} \leq
\|\phi\|_{L_p}+\|\phi_{xx}\|_{H_p^{-1}}$, (\ref{prelim-est1}) and
(\ref{prelim-est2}), we get:
\begin{eqnarray*}
\|u\|_{\bH_p^{1}} \leq \|u\|_{\bH_p^{0}}+ \|u_{xx}\|_{\bH_p^{-1}}
\leq  N ( \|f\|_{\bH_{p}^{-1}}+\|g\|_{\bL_{H}^{1,p}(L_p,l_2)}).
\end{eqnarray*}

\noindent We conclude that $u \in \bH_p^1$  and $u_{xx} \in
\bH_{p}^{-1}$, and hence $u \in \cH_{p,H}^1$.
Since ${\bf D}u=\Delta u+f$, we also infer that $\|u\|_{\cH_{p,H}^1}
\leq   N (\|f\|_{\bH_{p}^{-1}}+ \|g\|_{\bL_{H}^{1,p}(L_p,l_2)})$.

{\em Case 2.} The case of arbitrary $g=(g^k)_k \in
\widetilde{\bL}_{H}^{1,p}(L_p,l_2)$ follows as in the proof of
Theorem 4.2 of \cite{krylov99}, using an approximation argument.
This is based on the validity of the result in {\em Case 1} and the
completeness of the spaces $\bH_{p}^{n-2}$,
$\widetilde{\bL}_{H}^{1,p}(H_p^{n-1},l_2)$ and $\cH_{p,H}^n$
(Theorem \ref{krylov-theorem3-7}.(b)) $\Box$

\vspace{3mm}

Recall that, if $V$ is a Banach space and $\sigma \in (0,1)$, the
H\"older space $C^{\sigma}([0,T],V)$ is defined as the class of all
continuous functions $u :[0,T] \to V$ with
$$\|u\|_{C^{\sigma}([0,T],V)}:=\sup_{t \in [0,T]}\|u(t)\|_{V}+
\sup_{0\leq s<t \leq T}
\frac{\|u(t)-u(s)\|_{V}}{(t-s)^{\sigma}}<\infty.$$

Our final result is an embedding theorem for the space
$\cH_{p,H}^n$, similar to Theorem 7.2 of \cite{krylov99}.

\begin{theorem}
Let $p>2$, $n \in \bR$ and $1/2 \geq \beta>\alpha>1/p$. If $u \in
\cH_{p,H}^{n}$ then $u \in C^{\alpha-1/p}([0,T],H_p^{n-2\beta})$
a.s. and
$$E\|u(t,\cdot)-u(s,\cdot)\|_{H_p^{n-2\beta}}^p \leq  N(d,\beta,p,T)
(t-s)^{\beta p-1} \|u\|_{\cH_{p,H}^{n}}^p, \quad \forall 0 \leq s<t
\leq T;
$$
$$E\|u\|_{C^{\alpha-1/p}([0,T],H_p^{n-2\beta})}^{p} \leq
N(d,\beta,\alpha,p,T) \|u\|_{\cH_{p,H}^{n}}^p.$$
\end{theorem}

\noindent {\bf Proof:} 
We define $f={\bf D}u-\Delta u$, $g={\bf S}u$ and $u_0=u(0,\cdot)$.
Then $u$ satisfies the equation $dv=(\Delta v+f)dt +\sum_k g^k
\delta \beta_t^k$, with initial condition $v(0,\cdot)=u_0$. By
Theorem \ref{main-theorem}, this equation has a unique solution $v
\in \cH_{p,H}^n$. It follows that $u(t,\cdot)=v(t,\cdot)$ for all $t
\in [0,T]$, and it suffices to prove the theorem for $v$ in place of
$u$. By Proposition \ref{(1-Delta)-isometry}, without loss of
generality, we take $n=2\beta$. The theorem will be proved once we
show that
\begin{equation}
\label{embedding-ineq1-u} E\|u(t, \cdot)-u(s,\cdot)\|_{L_p}^p \leq
N(t-s)^{\alpha p-1} \{\|f\|_{\bH_{p}^{n-2}}^p+
\|g\|_{\bL_{H}^{1,p}(H_p^{n-1},l_2)}^p+E\|u_0\|_{H_p^{n-2/p}}^{p}\}
\end{equation}
\begin{equation}
\label{embedding-ineq2-u} E\sup_{0 \leq s<t \leq T}
\frac{\|u(t,\cdot)-u(s,\cdot)\|_{L_p}^p}{(t-s)^{\alpha p-1}} \leq N
\{\|f\|_{\bH_{p}^{n-2}}^p+
\|g\|_{\bL_{H}^{1,p}(H_p^{n-1},l_2)}^p+E\|u_0\|_{H_p^{n-2/p}}^{p}\}.
\end{equation}

Using an approximation argument and Theorem \ref{krylov-theorem3-7},
it is enough to assume that $u_0(\cdot)=1_{A_0}\phi(\cdot)$ with
$A_0 \in \cF, \phi \in C_0^{\infty}$,
\begin{equation}
\label{special-form-f-g} f(t,\cdot)=\sum_{i=1}^{m}
\sum_{j=1}^{m'}1_{A_j}1_{(t_{i-1},t_{i}]}(t)f_{ij}(\cdot) \quad
\mbox{and} \quad g^k(t,\cdot)=\sum_{i=1}^{m_k}F_i^k
1_{(t_{i-1}^{k},t_i^k]}(t) g_i^k(\cdot)
\end{equation}
where $A_j \in \cF$, $0\leq t_1<\ldots<t_m \leq T$ (non-random),
$f_{ij} \in C_0^{\infty}$, $F_i^k \in \cS_{\beta^k}$, $0\leq
t_1^k<\ldots<t_{m_k}^k \leq T$ (non-random), $g_i^k \in
C_0^{\infty}$, and $g_{i}^{k}=0$ for $k>K$.

Clearly, $u_0 \in L_p(\Omega,\cF, H_p^{2-2/p})$, $f \in \bH_p^{0}$
and $g \in \bL_{H}^{1,p}(H_p^1,l_2)$. By Theorem \ref{main-theorem},
it follows that $u \in \cH_{p,H}^2$. By Theorem
\ref{krylov-theorem3-7}.(a), $u \in C([0,T],L_p)$ a.s.

Let $u_1(t,x)=T_tu_0(x)+\int_0^t T_{t-s}f(s,\cdot)(x)ds$ and
$u_2(t,x)=u(t,x)-u_1(t,x)$. Relations (\ref{embedding-ineq1-u}) and
(\ref{embedding-ineq2-u}) for $u_1$ follow as in the proof of
Theorem 7.2 of \cite{krylov99}.

Hence, it suffices to prove (\ref{embedding-ineq1-u}) and
(\ref{embedding-ineq2-u}) for $u_2$. Using (\ref{formula-for-u2}),
it follows that
$$u_2(r+\gamma,x)-u_2(r,x)=(T_{\gamma}-1)u_2(r,\cdot)(x)+\sum_{k=1}^{\infty}
\int_{r}^{r+\gamma} T_{r+\gamma-\rho}g^{k}(\rho,\cdot)(x) \delta
\beta_{\rho}^{k}$$ and hence $E\|u_2(r+\gamma,
\cdot)-u_2(r,\cdot)\|_{L_p}^p \leq N (A_2(r,\gamma)+B_2(r,\gamma))$,
where
\begin{eqnarray*}
A_2(r,\gamma) & :=& E \int_{\bR^d} |(T_{\gamma}-1)u_2(r,\cdot)(x)|^p
dx \\
B_2(r,\gamma)  &:= & E \int_{\bR^d}
\left|\sum_{k=1}^{\infty}\int_{r}^{r+\gamma}
T_{r+\gamma-\rho}g^k(\rho,\cdot)(x) \delta \beta_{\rho}^k
\right|^{p}dx.
\end{eqnarray*}

\noindent We now apply Lemma 7.4 of \cite{krylov99} to the
continuous function $u_2:[0,T] \to L_p$:
$$E\|u_2(t,\cdot)-u_2(s,\cdot)\|_{L_p}^p \leq N(t-s)^{\alpha p-1}
(I_{2}(t,s)+J_2(t,s))$$
$$E \sup_{0 \leq s<t \leq T} \frac{\|u_2(t,\cdot)-u_2(s,\cdot)\|_{L_p}^{p}}
{(t-s)^{\alpha p-1}} \leq N (I_{2}(t,s)+J_2(t,s)),$$ with
$$I_2(t,s)=\int_0^{t-s} \frac{d\gamma}{\gamma^{1+\alpha p}}
\int_{s}^{t-\gamma} A_2(r,\gamma)dr, \ J_2(t,s)=\int_0^{t-s}
\frac{d\gamma}{\gamma^{1+\alpha p}} \int_{s}^{t-\gamma}
B_2(r,\gamma)dr.$$ The term $I_2(t,s)$ is estimated as in
\cite{krylov99}, using Theorem \ref{main-theorem}:
\begin{equation}
\label{estimate-I2} I_2(t,s) \leq N (t-s)^{(\beta-\alpha)p}
\|g\|_{\bL_{H}^{1,p}(H_p^{n-1},l_2)}^{p}.
\end{equation}

It remains to estimate $J_2(t,s)$. Using Theorem
\ref{max-inequality-p}, we have:
\begin{eqnarray}
\nonumber B_2(r,\gamma) & \leq & N \left\{\int_{\bR^d} E
\left(\int_{r}^{r+\gamma} |T_{r+\gamma-\rho}g(\rho,
\cdot)(x)|_{l_2}^{2} d\rho \right)^{p/2}dx +
\right. \\
\nonumber & & \left. \int_{\bR^{d}} E \left[\int_{r}^{r+\gamma}
\left( \int_{0}^{T} |D_{\theta}[T_{r+\gamma-\rho}
g(\rho,\cdot)(x)]|_{l_2}^{1/H} d\theta \right)^{2H}d\rho
\right]^{p/2}dx \right\} \\
\label{estimate-B2} & := & N(B_{2}'(r,\gamma)+B_{2}''(r,\gamma)).
\end{eqnarray}

\noindent The term $B_{2}'(r,\gamma)$ is treated as in
\cite{krylov99}:
\begin{equation}
\label{estimate-B2'} B_{2}'(r,\gamma) \leq N \gamma^{\beta p-1} E
\int_{0}^{\gamma} \|g(r+\rho, \cdot)\|_{H_p^{n+1}(l_2)}^{p} d\rho.
\end{equation}

\noindent For the term $B_{2}''(r,\gamma)$, we use
(\ref{D-commutes-with-T}), H\"older's inequality with $q=p/(p-2)$,
Minkowski's inequality, and Lemma 7.3 of \cite{krylov99}:
\begin{eqnarray}
\nonumber \lefteqn{ B_{2}''(r,\gamma) =E \left[\int_{r}^{r+\gamma}
\left( \int_{0}^{T} |T_{r+\gamma-\rho}
[D_{\theta}g(\rho,\cdot)](x)|_{l_2}^{1/H} d\theta \right)^{2H}d\rho
\right]^{p/2}dx } \\
\nonumber & & = E \left[\int_{0}^{\gamma} \rho^{2\beta-1}
\rho^{1-2\beta} \left(\int_0^T |T_{\rho}
[D_{\theta}g(r+\gamma-\rho,\cdot)](x)|_{l_2}^{1/H}d\theta
\right)^{2H} d\rho \right]^{p/2}dx \\
\nonumber & & \leq N \gamma^{\beta p-1} E \int_0^{\gamma}
\rho^{(1-2\beta)p/2} \int_{\bR^d} \left(\int_0^T |T_{\rho}
[D_{\theta}
g(r+\gamma-\rho,\cdot)](x)|_{l_2}^{1/H}d\theta\right)^{pH}dx d\rho \\
\nonumber & & \leq N \gamma^{\beta p-1} E \int_0^{\gamma}
\rho^{(1-2\beta)p/2} \left[\int_0^T \left(\int_{\bR^d} |T_{\rho}
[D_{\theta} g(r+\gamma-\rho,\cdot)](x)|_{l_2}^{p} dx
\right)^{1/(pH)}d\theta
\right]^{pH}d\rho \\
\nonumber & & \leq N \gamma^{\beta p-1} E
\int_0^{\gamma}\rho^{(1-2\beta)p/2}
\left(\frac{e^{\rho}}{\rho^{1/2-\beta}} \right)^p\left(\int_0^{T} \|
D_{\theta} g(r+\gamma-\rho,\cdot)\|_{H_{p}^{n-1}(l_2)}^{1/H}d\theta
\right)^{pH}d\rho \\
\label{estimate-B2''} & & = N \gamma^{\beta p-1} E \int_0^{\gamma}
\left(\int_0^{T} \| D_{\theta}
g(r+\gamma-\rho,\cdot)\|_{H_{p}^{n-1}(l_2)}^{1/H}d\theta
\right)^{pH}d\rho.
\end{eqnarray}

Using (\ref{estimate-B2}), (\ref{estimate-B2'}) and
(\ref{estimate-B2''}), we obtain:
\begin{eqnarray}
\nonumber \lefteqn{J_2(t,s) \leq  N \left\{ E \int_0^{t-s}
\frac{1}{\gamma^{2+(\alpha-\beta) p}}  \int_{s}^{t-\gamma}
\int_0^{\gamma}
\|g(r+\rho,\cdot)\|_{H_p^{n-1}(l_2)}^{p} dr d\rho d\gamma + \right. } \\
\nonumber & & \left. E \int_0^{t-s}
\frac{1}{\gamma^{2+(\alpha-\beta) p}} \int_{s}^{t-\gamma}
\int_0^{\gamma} \left(\int_0^T
\|D_{\theta}g(r+\rho,\cdot)\|_{H_p^{n-1}(l_2)}^{1/H}d\theta
\right)^{pH}
dr d\rho d\gamma \right\} \\
\nonumber &  & \leq N (t-s)^{(\beta-\alpha)p} \left\{ E \int_0^{t}
\|g(r,\cdot)\|_{H_p^{n-1}(l_2)}^{p} dr + E \int_0^{t} \left(\int_0^T
\|D_{\theta}g(r,\cdot)\|_{H_p^{n-1}(l_2)}^{1/H}d\theta \right)^{pH}
dr \right\} \\
\label{estimate-J2} & & \leq N (t-s)^{(\beta-\alpha)p}
\|g\|_{\bL_{H}^{1,p}(H_p^{n-1},l_2)}^{p}.
\end{eqnarray}

Relations (\ref{embedding-ineq1-u}) and (\ref{embedding-ineq2-u})
for $u_2$ follow from (\ref{estimate-I2}) and (\ref{estimate-J2}).
$\Box$

\appendix

\section{A Banach-space generalization of Littlewood-Paley inequality}

Let $V$ be an arbitrary Hilbert space. For any $f \in
L_p(V)=L_p(\bR^d,V)$, $p \geq 1$, we let $$T_t f (x):=\int_{\bR^d}
f(x-y)G_t(y)dy,$$ where $G_t(x)=(4\pi t)^{-d/2}\exp\{-|x|^2/(4t)\},
t>0, x \in \bR^d$ is the heat kernel.

First, notice that:
\begin{equation}
\label{ineq-Tt-g-norm-p-V} \|T_t f \|_{L_p(V)} \leq \|f\|_{L_p(V)}.
\end{equation}

To see this, note that $|T_t f (x)|_{V} \leq \int_{\bR^d}
|f(x-y)|_{V} G_t(y)dy$ for any $x \in \bR^d$. Using Minkowski's
inequality for integrals, we have:
\begin{eqnarray*}
 \|T_t f\|_{L_{p}(V)} & \leq &
\left[\int_{\bR^d} \left(\int_{\bR^d}
|f(x-y)|_{V}G_t(y)dy \right)^{p}dx\right]^{1/p}  \\
& \leq & \int_{\bR^d}G_t(y) \left( \int_{\bR^d}
|f(x-y)|_{V}^{p}dx\right)^{1/p}dy =\|f\|_{L_p(V)} \|G_t\|_{L_1}.
\end{eqnarray*}

The following result is a generalization of the Littlewood-Paley
inequality, due to \cite{krylov94} (see Theorem 1.1 of
\cite{krylov94}, and \cite{krylov06}).

\begin{theorem}\label{theorem1-1-krylov} Let $p \in [2,\infty)$
and $f \in C_{0}^{\infty}((a,b) \times \bR^d,V)$, where $-\infty
\leq a<b \leq \infty$. Then
$$\int_{\bR^d} \int_a^b \left[ \int_a^t | \triangledown
T_{t-s}f(s, \cdot)(x)|_{V}^{2} ds \right]^{p/2} dt dx \leq N
\int_{\bR^d} \int_a^b |f(t,x)|_{V}^{p}dt dx,$$ where $N$ is a
constant depending only on $d$ and $p$.
\end{theorem}

In the present article, we need the following generalization of
Theorem \ref{theorem1-1-krylov} to the case of $U$-valued functions,
where $U=L_{1/H}((\alpha,\beta),V)$ is a Banach space.

\begin{theorem}
\label{gener-Th1-1-Krylov} Let $p \in [2, \infty)$ and $f \in
C_0^{\infty}((a,b) \times \bR^d,U)$, where  $-\infty \leq a<b\leq
\infty$ and $U=L_{1/H}((\alpha,\beta),V)$, with $-\infty \leq
\alpha<\beta \leq \infty$ and $1/2<H<1$. Then
\begin{eqnarray}
\nonumber \lefteqn{\int_{\bR^d} \int_a^b \left[ \int_a^t \left(
\int_{\alpha}^{\beta} |\triangledown T_{t-s} f(s, \cdot,
\theta)(x)|_{V}^{1/H}d\theta \right)^{2H} ds \right]^{p/2}
dt dx \leq } \\
\label{LHS-gener-Th1-1-Krylov} & & N \int_a^b \left[
\int_{\alpha}^{\beta} \left( \int_{\bR^d}|f(t,x,\theta)|_{V}^{p} dx
\right)^{1/pH} d\theta \right]^{pH}dt,
\end{eqnarray}
where $N$ is a constant depending only on $d$ and $p$.
\end{theorem}

The remaining part of this section is dedicated to the proof of
Theorem \ref{gener-Th1-1-Krylov}. We follow the lines of the proof
of Theorem 16.1 of \cite{krylov06}. It is enough to assume that
$a=-\infty$ and $b=\infty$. We first treat the case $p=2$.

\begin{lemma}
\label{BK-lema1} Relation (\ref{LHS-gener-Th1-1-Krylov}) holds for
$p=2$.
\end{lemma}

\noindent {\bf Proof:}
Due to Minkowski's inequality, the left-hand side of
(\ref{LHS-gener-Th1-1-Krylov}) is smaller than
$$\int_{-\infty}^{\infty} \left[\int_{\alpha}^{\beta} \left(
\int_s^{\infty} \int_{\bR^d} |\triangledown T_{t-s} f(s,\cdot,
\theta)(x)|_{V}^2dx dt\right)^{1/(2H)}d\theta \right]^{2H}ds.$$
Using the Fourier transform, the inner integral equals
\begin{eqnarray*}
\int_s^{\infty} \int_{\bR^d}|\xi|^2 e^{-(t-s)|\xi|^2} |{\cF}f(s,
\xi,\theta)|_{V}^2 d\xi dt&=&\int_{\bR^d}|\cF f(s,\xi,\theta)|_{V}^2
|\xi|^2
\left(\int_s^{\infty}e^{-(t-s)|\xi|^2}dt \right) d\xi \\
&=& \int_{\bR^d}|\cF f(s,\xi,\theta)|_{V}^2 d\xi,
\end{eqnarray*}
which proves (\ref{LHS-gener-Th1-1-Krylov}) for $p=2$. $\Box$

\vspace{3mm}

Assume now that $p>2$. Note that $\triangledown T_t
h(x)=t^{-1/2}\Psi_t h(x)$, where $\Psi_t
h(x)=t^{-d/2}\phi(x/\sqrt{t})*h(x)$ and $\phi(x)=-(4\pi)^{-d/2} x
e^{-|x|^2/4}$. Set
\begin{eqnarray*}
u(t,x)=\cG f(t,x)&=&\left[\int_{-\infty}^{t}
\left(\int_{\alpha}^{\beta} |\Psi_{t-s}f(s,\cdot,
\theta)(x)|^{1/H}d\theta \right)^{2H} \frac{1}{t-s}ds\right]^{1/2}.
\\
&=& \left(\int_{-\infty}^{t} |\Psi_{t-s}f(s,\cdot, *)(x)|_{U}^2 \
\frac{1}{t-s}ds\right)^{1/2}.
\end{eqnarray*}

We want to prove that:
\begin{equation}
\label{Krylov-Banach} \int_{\bR^d} \int_{-\infty}^{\infty}
|u(t,x)|^p dt dx \leq N \int_{-\infty}^{\infty} \left[ \left(
\int_{\bR^d} |f|_{V}^p(t,x,\theta)dx
\right)^{1/(pH)}d\theta\right]^{pH}dt.
\end{equation}


Recall that the {\bf maximal function} of $g:\bR^d \to \bR$ is
defined by:
$$\bM_x g(x)=\sup_{r>0} \frac{1}{B_r} \int_{B_r(x)}|g(y)|dy,$$
where $B_r(x)=\{y; |y-x| <r\}$ and $B_r=B_r(0)$. If $h:\bR^{d+1} \to
\bR$, we define $\bM_{x} h(t,x)=\bM_x h(t,\cdot)(x)$. Let
$Q_0=[-4,0] \times [-1,1]^d$.

\begin{lemma}
\label{BK-lema2} Assume that $f(t,x,\theta)=0$ for $(t,x) \not \in
(-12,12) \times B_{3d}$. Then for any $(t,x) \in Q_0$
\begin{equation}
\label{BK-lema23-rel} \int_{Q_0}|u(s,y)|^2 dsdy \leq N \bM_{t} \|
\bM_{x} |f|_{V}^2 (t,x,*) \|_{U_0},
\end{equation}
where $U_0=L_{1/(2H)}((\alpha,\beta))$ and $N$ depends only on $d$.
\end{lemma}

\noindent {\bf Proof:} Using Lemma \ref{BK-lema1}, the left-hand
side of (\ref{BK-lema23-rel}) is smaller than:
$$N\int_{-\infty}^0 \left[\int_{\alpha}^{\beta}\left(\int_{\bR^d}
|f|_{V}^2(s,y,\theta) dy \right)^{1/(2H)} d\theta\right]^{2H}ds \leq
N  \int_{-12}^{0}
\left[\int_{\alpha}^{\beta}(\bM_{x}|f|_{V}^2(s,x,\theta))^{1/(2H)}
d\theta\right]^{2H}ds$$
$$=N  \int_{-12}^{0}
\|\bM_{x}|f|_{V}^2(s,x,*) \|_{U_0}ds \leq N \bM_{t} \| \bM_{x}
|f|_{V}^2 (t,x,*) \|_{U_0}.$$ $\Box$

\begin{lemma}
\label{BK-lema3} Assume that $f(t,x,\theta)=0$ for $t \not \in
(-12,12)$. Then (\ref{BK-lema23-rel}) holds for any $(t,x) \in Q_0$.
\end{lemma}

\noindent {\bf Proof:} Let $\zeta \in C_0^{\infty}(\bR^d)$ be such
that $\zeta=1$ in $B_{2d}$, $\zeta=0$ outside $B_{3d}$, and
$\zeta(x) \in (0,1)$ for $x \in B_{3d} \verb2\2 B_{2d}$. Let
$\alpha=\zeta f$ and $\beta =(1-\zeta)f$. Then
\begin{eqnarray*}
\cG f(t,x) &=& \cG (\alpha+\beta)(t,x) = \left(\int_{-\infty}^{t}
|\Psi_{t-s}(\alpha+\beta)(s,\cdot,*)(x)|_{U}^2
\frac{1}{t-s}ds\right)^{1/2} \\
& \leq & \cG \alpha(t,x)+\cG \beta(t,x),
\end{eqnarray*}
using Minkowski's inequality in $L_{2}(\bR,U)$, which in turn relies
on Minkowski's inequality in the Banach space $U$. Since $\alpha$
satisfies the conditions of Lemma \ref{BK-lema2} and $|\alpha|_V
\leq |f|_{V}$, for any $(t,x) \in Q_0$
$$\int_{Q_0} |\cG\alpha(s,y)|^2dsdy \leq
N \bM_{t} \|\bM_{x} |\alpha|_{V}^{2} (t,x,*) \|_{U_0} \leq N \bM_{t}
\|\bM_{x} |f|_{V}^{2} (t,x,*) \|_{U_0}.$$

Therefore, it suffices to prove that (\ref{BK-lema23-rel}) holds for
any function $f$ such that $f(t,x,\theta)=0$ if $t \not \in
(-12,12)$ or $x \in B_{2d}$ (in particular for $\beta$). This
follows as in the proof of Lemma 16.5 of \cite{krylov06}, using
Minkowski's inequality for integrals. $\Box$

\begin{lemma}
\label{BK-lema4} Assume that $f(t,x,\theta)=0$ for $t \geq -8$. Then
for any $(t,x) \in Q_0$
$$\int_{Q_0}|u(s,y)-u(t,x)|^2 dsdy \leq N \bM_{t}
\|\bM_{x} |f|_{V}^{2} (t,x,*) \|_{U_0}.$$
\end{lemma}

\noindent {\bf Proof:} The argument is similar to the one used in
the proof of Lemma 16.6 of \cite{krylov06}, with some minor
modifications (as above). $\Box$

\vspace{3mm}

We introduce now the filtration $\bQ_n, n \in \bZ$ of partitions
$\bQ_n=\{Q_n(i_0,i_1, \ldots,i_d); \linebreak i_0,i_1\ldots,i_d \in
\bZ\}$ of $\bR^{d+1}$, as in \cite{krylov06}. For any $x \in \bR^d$
and $n \in \bZ$, we denote by $Q_n(x)$ the unique $Q \in \bQ_n$
containing $x$. The {\bf sharp function} of $g \in L_{1,
loc}(\bR^d)$ is defined by:
$$g^{\#}(x)=\sup_{n \in \bZ} \frac{1}{|Q_n(x)|}
\int_{Q_n(x)}|g(y)-g_{|n}(x)|dy,$$ where $g_{|n}(x)=|Q_n(x)|^{-1}
\int_{Q_n(x)}g(y)dy$. If $p \in (1,\infty)$, then by the
Fefferman-Stein theorem, for any $g \in L_{p}(\bR^d)$,
$\|g\|_{L_p(\bR^{d})} \leq N \|g^{\#}\|_{L_p(\bR^d)}$.

\begin{lemma}
\label{BK-lema5} Let $f \in C_0^{\infty}(\bR^{d+1},U)$ be arbitrary.
For any $(t,x) \in \bR^{d+1}$,
$$(\cG f)^{\#}(t,x) \leq N (\bM_{t}
\|\bM_{x} |f|_{V}^{2} (t,x,*) \|_{U_0})^{1/2}.$$
\end{lemma}

\noindent {\bf Proof:} The argument is based on Lemma \ref{BK-lema3}
and Lemma \ref{BK-lema4}, and is similar to the one used for proving
relation (16.20) of \cite{krylov06}. $\Box$

\vspace{3mm}

\noindent {\bf Proof of Theorem \ref{gener-Th1-1-Krylov}:} Assume
that $p>2$. We use the Fefferman-Stein theorem, Lemma
\ref{BK-lema5}, the boundedness of the operators $\bM_{t}$ and
$\bM_x$ ($p>2$), and Minkowski's inequality for integrals ($pH>1$):
\begin{eqnarray*}
\|u\|_{L_{p}(\bR^{d+1})}^p & \leq & N \|(\cG f)^{\#}
\|_{L_p(\bR^{d+1})}^p \leq N \int_{\bR^d} \int_{\bR}(\bM_{t}
\|\bM_{x} |f|_{V}^{2} (t,x,*) \|_{U_0})^{p/2}dt dx \\
&=& N \int_{\bR^d} \| \bM_{t} \|\bM_{x} |f|_{V}^{2} (t,x,*)
\|_{U_0}\|_{L_{p/2}(\bR)}^{p/2}dx \\
& \leq & N \int_{\bR^d} \int_{\bR} \|\bM_{x} |f|_{V}^{2} (t,x,*)
\|_{U_0}dt dx \\
&=& N \int_{\bR} \int_{\bR^d} \left[\int_{\alpha}^{\beta}(\bM_{x}
|f|_{V}^{2} (t,x,\theta))^{1/(2H)} d\theta \right]^{pH} dx dt \\
&\leq & N \int_{\bR} \left[\int_{\alpha}^{\beta} \left(\int_{\bR^d}
(\bM_{x} |f|_{V}^{2} (t,x,\theta))^{p/2} dx \right)^{1/(pH)}d\theta
\right]^{pH} dt \\
&\leq & N \int_{\bR} \left[\int_{\alpha}^{\beta} \left(\int_{\bR^d}
|f|_{V}^{p} (t,x,\theta)) dx \right)^{1/(pH)}d\theta \right]^{pH}
dt,
\end{eqnarray*}
i.e. (\ref{Krylov-Banach}) holds. $\Box$

\noindent \
University of Ottawa \\
Department of Mathematics and Statistics \\
585 King Edward Avenue
Ottawa, ON, K1N 6N5, Canada \\
E-mail address: rbalan@uottawa.ca \\
URL: http://aix1.uottawa.ca/\~ \ rbalan

\end{document}